\documentclass[reqno,a4paper]{amsart}

\usepackage{graphicx}
\usepackage{color} 
\usepackage[usenames,dvipsnames]{xcolor}
\usepackage{transparent}
\usepackage[utf8x]{inputenc}  
\usepackage{dsfont}
\usepackage{gensymb}
\usepackage[active]{srcltx}\usepackage[colorlinks,pdfpagelabels,pdfstartview=FitH,bookmarksopen=true,bookmarksnumbered=true,linkcolor=blue,plainpages=false,hypertexnames=false,citecolor=red]{hyperref}
\usepackage{mathrsfs}
\usepackage{verbatim, amssymb,amsmath, amsfonts, amsthm}
\usepackage{latexsym}
\usepackage[dvips]{epsfig}
\usepackage{float}

\allowdisplaybreaks 

\newtheorem{teorema}{Theorem}[section]
\newtheorem{proposizione}[teorema]{Proposition}
\newtheorem{lemma}[teorema]{Lemma}
\newtheorem{corollario}[teorema]{Corollary}
\newtheorem{osservazione}[teorema]{Remark}
\newtheorem{definizione}[teorema]{Definition}



\newcommand{\R}{\mathbb{R}}
\renewcommand{\to}{\rightarrow}

\def\sideremark#1{\ifvmode\leavevmode\fi\vadjust{\vbox to0pt{\vss
			\hbox to 0pt{\hskip\hsize\hskip1em
				\vbox{\hsize2.1cm\tiny\raggedright\pretolerance10000
					\noindent #1\hfill}\hss}\vbox to15pt{\vfil}\vss}}}%

\newcommand{\virg}[1]{\textquotedblleft#1\textquotedblright}

\begin{document}
	\numberwithin{equation}{section}
	\parindent=0pt
	\hfuzz=2pt
	\frenchspacing

	\title[]{Stability and asymptotic behaviour of one-dimensional solutions in cylinders}
	
	\author[]{Francesca De Marchis, Lisa Mazzuoli, Filomena Pacella}
	
	\address{Francesca De Marchis, University of Roma {\em Sapienza}, P.le Aldo Moro 5, 00185 Roma, Italy -- demarchis@mat.uniroma1.it}
	\address{Lisa Mazzuoli, University of Roma {\em Sapienza}, P.le Aldo Moro 5, 00185 Roma, Italy -- lisa.mazzuoli@uniroma1.it}
	\address{Filomena Pacella, University of Roma {\em Sapienza}, P.le Aldo Moro 5, 00185 Roma, Italy -- filomena.pacella@uniroma1.it}
	
	
	
	\thanks{Research partially supported by PRIN Project 2022AKNSE4 \emph{Variational and Analytical aspects of Geometric PDE} and by Gruppo Nazionale per l’Analisi Matematica, la Probabilità e le loro
		Applicazioni (GNAMPA) of the Istituto Nazionale di Alta Matematica (INdAM). In particular the first author is supported by INdAM-GNAMPA Project $\#$E5324001950001$\#$.}
	
	\begin{abstract} 
		We consider positive one-dimensional solutions of a Lane-Emden relative Dirichlet problem in a cylinder and study their stability/instability properties as the energy varies with respect to domain perturbations. This depends on the exponent $p >1$ of the nonlinearity and we obtain results for $p$ close to 1 and for $p$ large. This is achieved by a careful asymptotic analysis of the one-dimensional solution as $p \to 1$ or $p \to \infty$, which is of independent interest. It allows to detect the limit profile and other qualitative properties of these solutions.
	\end{abstract}
	
	\maketitle
	\section{Introduction}
	In this paper we study one-dimensional 
	solutions of a class of semilinear elliptic problems in $\mathbb{R}^N, N \geq 2$. By this we mean functions $u$ which depends only on one cartesian variable, say the $x_N$-variable, i.e. if $x=(x_1,...,x_N)$ is a point in $\mathbb{R}^N$, then 
	\begin{equation*}
		u(x)=u(x_1,...x_N)=v(x_N)
	\end{equation*}
	for some function $v$ of one variable. \\
	These functions need to be defined in specific domains, namely cylinders, either bounder or unbounded. Indeed a cylinder is defined as
	\begin{equation*}
		\Sigma_\omega=\omega \times I
	\end{equation*}
	where $\omega$ is a bounded domain in $\mathbb{R}^{N-1}$ and $I$ is an interval in $\mathbb{R}$. Thus if $x=(x',x_N),$ $x' \in \mathbb{R}^{N-1}$ the $x_N$-variable represents a distinguished direction which allows to define the one-dimensional functions in $\Sigma_\omega$. Note that these functions have flat level sets. They arise in several PDE problems as for example in connection with the famous De Giorgi conjecture (see \cite{GhossGuis} and the references therein).
	However, in spite of their simple form, we believe that their properties in the study of PDEs have not been much  investigated.\\
	In this paper we consider unbounded semi-cylinders $\Sigma_\omega$ with a smooth bounded cross-section $\omega$, i.e.
	\begin{equation*}
		C_\omega=\omega \times (0, +\infty)
	\end{equation*}
	where $\omega$ is a bounded smooth domain in $\mathbb{R}^{N-1}$, $N\geq 2$. If $\Omega$ is a bounded domain in $C_\omega$, we denote by $\Gamma_\Omega$ its relative (or free) boundary in $C_\omega$:
	\begin{equation*}
		\Gamma_{\Omega}= \partial \Omega \cap C_\omega
	\end{equation*}
	and by $\Gamma_{1,\Omega}$ the part of $\partial \Omega$ on $\partial C_\omega$, i.e.
	\begin{equation*}
		\Gamma_{1,\Omega}=\partial \Omega \cap \partial C_\omega.
	\end{equation*}
	We assume that $\Gamma_\Omega$ is smooth.\\
	We focus on the relative Dirichlet Lane-Emden problem which is the following
	\begin{equation} \label{problem_Omega}
		\begin{cases}
			-\Delta u=u^p \quad &\mbox{ in }\Omega\\
			u>0  &\mbox{in } \Omega\\
			u=0 &\mbox{on } \Gamma_\Omega\\
			\frac{\partial u}{\partial \nu}=0 & \mbox{on } \Gamma_{1,\Omega}
		\end{cases}
	\end{equation}
	where $p>1$.
	By standard variational methods it is easy to see that a positive weak solution of $\eqref{problem_Omega}$ exists in the Sobolev space $H^1_0(\Omega \cup \Gamma_{1,\Omega})$, which is the space of functions in $H^1(\Omega)$ whose trace vanishes on the relative boundary $\Gamma_\Omega$, at least for some values of $p$.\\
	Let $u_\Omega$ be a positive solution of $\eqref{problem_Omega}$, then its energy is defined as:
	\begin{equation}\label{energy_functional}
		J(u_\Omega)=\frac{1}{2}\int_{\Omega} |\nabla u_\Omega |^2\,dx -\frac{1}{p+1}\int_{\Omega}u_{\Omega}^{p+1}\, dx.
	\end{equation}
	An interesting question is to understand how the energy $J(u_\Omega)$ behaves with respect to variations of $\Omega$ which preserve its measure.\\
	Roughly speaking, we consider the functional 
	\begin{equation}\label{T_Omega}
		T(\Omega)=J(u_\Omega)
	\end{equation}
	and study it under small perturbations of $\Omega$. To well define the functional $T$ a local uniqueness result of the positive solution $u_\Omega$  is needed and this is ensured by considering a nondegenerate solution $u_\Omega$ (see Section $\ref{section:stability}$ and \cite[Proposition 2.1]{AfonsoPacella}).\\
	In particular, domains which are local minima of T are of special interest. This can be considered a local shape-optimization problem. Let us consider in $C_\omega$ domains $\Omega_\varphi$ which are the hypographs of functions $e^\varphi$, $\varphi \in C^2(\bar{\omega})$, i.e.
	\begin{equation*}
		\Omega_\varphi= \lbrace x=(x',x_N) \in C_\omega: \quad x_N< e^{\varphi(x')} \rbrace,
	\end{equation*}
	for a function $\varphi \in C^2(\bar{\omega})$.\\
	Note that the relative boundary of $\Omega_\varphi$ is just the cartesian graph of $e^\varphi$ and we denote it simply by $\Gamma_{\varphi}$.\\\\
	A particular simple domain is obtained by taking the hypograph of a constant function. It will be defined as:
	\begin{equation} \label{Omega}
		\Omega_{\varphi_L}=\lbrace (x',x_N), x' \in \omega, \,\, 0<x_N<L\rbrace,
	\end{equation}
	hence $\varphi_L=\log L$, $L>0$.\\
	Obviously $\Omega_{\varphi_L}$ is just the bounded cylinder $\omega \times (0,L)$ with base $\omega$ and height $L$.\\
	In this domain there exists a one-dimensional solution of $\eqref{problem_Omega}$ which is easily obtained by extending to $\Omega_{\varphi_L}$ a positive solution $u_{p,L}$ of the following ODE problem:
	\begin{equation}\label{up_problem}
		\begin{cases}
			-u''=u^p \quad &\mbox{in } (0,L)\\
			u'(0)=u(L)=0
		\end{cases}
	\end{equation}
	i.e. with an abuse of notation we define 
	\begin{equation}\label{upLN}
		u_{p,L}(x',x_N):=u_{p,L}(x_N), \quad \forall\, x' \in \omega.
	\end{equation}
	In view of the simple geometry of $\Omega_{\varphi_L}$ and of the corresponding one-dimensional solution $u_{p,L}$, we may ask wether $\Omega_{\varphi_L}$ is or not a good candidate to locally minimize the energy functional $T$. This question is addressed in \cite{AfonsoPacella} (see also \cite{Afonso}) by differentiating $T$ with respect to variation of $\Omega$ which leave the measure invariant and studying the stability/instability of the pair $(\Omega_{\varphi_L}, u_{p,L})$ as a critical point of $T$ (see Section $\ref{section:stability}$ for the definition and the setting). Indeed, when $u_{p,L}$ is a nondegenerate solution of $\eqref{problem_Omega}$ (see Proposition $\ref{prop:nondegeneracy}$) we may consider the energy functional $T({\Omega_\varphi})$ for domains $\Omega_\varphi$ close to $\Omega_{\varphi_L}$.\\
	In the paper \cite{AfonsoPacella} some conditions for the stability of the pair $(\Omega_{\varphi_L}, v_L)$ ($v_L$ one-dimensional nondegenerate solution) were derived, for semilinear elliptic problems with a general nonlinearity $f(u)$, while for the case of the torsion problem, i.e. $f(u)=1$, sharp results both for instability and stability were proved (see \cite[Theorem 1.4]{AfonsoPacella}).
	{A related result is obtained in \cite{Afonso1} by showing a bifurcation result from 1-dimensional solutions of semilinear elliptic problems in bounded cylinders.}
	\\\\
	In the present paper we analyze the case of the Lane-Emden nonlinearity and for asymptotic values of the exponent $p$ we substantially improve the results of \cite{AfonsoPacella} getting also instability. \\
	Our results rely on a careful asymptotic analysis of the one-dimensional solution $u_{p,L}$, as $p\to \infty$ or as $p \to 1$, which is new and interesting in itself.\\
	As far as we know this is the first time that a qualitative analysis of the one-dimensional solutions is performed and we believe that it can be useful for other problems.\\\\
	For simplicity to state our asymptotic results we fix $L=1$ and denote $u_{p,1}$ simply by $u_p$. Moreover it is convenient to consider the even extension of $u_p$ to the interval $I=(-1,1)$, that, with an abuse of notation, we still denote by $u_p$.
	Clearly it is the unique solution of the Dirichlet problem 
	\begin{equation}\label{problema_in_I}
		\begin{cases}
			-u''=u^p \quad & \mbox{ in } I\\
			u>0 & \mbox{ in } I\\
			u(\pm 1)=0.
		\end{cases}
	\end{equation}
	By symmetry results $u_p$ is the only positive solution and it is decreasing in $[0,1]$, thus $\|u_p\|_\infty=u_p(0)$.\\
	Note that it is easy to pass from the case $L=1$ to the case of any other $L>0$ (see Remark \eqref{rem:equivalence}).\\
	We have
	\begin{teorema}[Asymptotic behaviour for $p \to +\infty$]\label{Teo_p_to_inf}
		Let $u_p$ be the positive solution of \eqref{problema_in_I} and let $\alpha_1(p)$ be the first eigenvalue of the linearized operator at $u_p$ with Dirichlet boundary conditions in $I$, defined as
		\begin{equation}\label{L_up}
			L_{u_p}=-\frac{d^2}{dt^2}-pu^{p-1}_{p}.
		\end{equation}
		Then we have
		\begin{enumerate}
			\item \begin{equation*}
				\lim_{p\to +\infty}\Vert u_p \Vert_{\infty}=1;
			\end{equation*}
			\item \begin{equation*}
				u_p\rightarrow 2G(\cdot,0) \qquad \mbox{ in } C^0([-1,1])\cap C^1_{loc}([-1,1]\setminus\lbrace 0 \rbrace),
			\end{equation*}
			where $G$ is the Green's function of $-\frac{d^2}{d t^2}$ with Dirichlet boundary conditions defined as
			\begin{equation}\label{Green}
				G(t,\tau)=\begin{cases}
					-\frac{(\tau+1)(t-1)}{2} \quad &\mbox{ if } t>\tau\\
					\frac{(t+1)(1-\tau)}{2}  \quad &\mbox{ if } t\leq \tau;
				\end{cases}
			\end{equation}
			\item 
			\begin{equation*}
				\Vert u_p \Vert_{\infty}^{p+1}=\frac{p}{2}(1+o_p(1)) \qquad \mbox{as $p\to+\infty$};
			\end{equation*}
			\item setting
			\[
			\mu_p:=(p\|u_p\|_\infty^{p-1})^{-\frac12}
			\]
			and
			\begin{equation}\label{utildep}
				\tilde{u}_p(s)=p\frac{u_p(\mu_ps)-u_p(0)}{\|u_p\|_\infty}, \qquad \text{for $s \in [-\tfrac{1}{\mu_p}, \tfrac{1}{\mu_p}]$,}
			\end{equation}
			then
			\begin{equation*}
				\tilde{u}_p \rightarrow W \quad \mbox{ in } C^1_{loc}(\mathbb{R})\qquad \mbox{as $p\to+\infty$}
			\end{equation*}
			where 
			\begin{equation}\label{W}
				W(s)=\log\frac{4e^{\sqrt{2}s}}{(1+e^{\sqrt{2}s})^2},
				\quad s\in \mathbb{R}.
			\end{equation}
			is the solution of the limit problem 
			\begin{equation}
				\begin{cases}
					-W''=e^W \quad \mbox{ in } \mathbb{R}\\
					W'(0)=0\\
					W(0)=0.
				\end{cases}
				\label{W2}
			\end{equation}
			\item
			\begin{equation*}
				\alpha_1(p)=-{\frac{1}{4}}p^2(1+o_p(1))\qquad \mbox{as $p\to+\infty$}.
			\end{equation*}
		\end{enumerate}
	\end{teorema}
	\begin{teorema}[Asymptotic behaviour for $p \to 1$]\label{Teo_p_to_1}
		Let $u_p$ be the positive solution of \eqref{problema_in_I}, let $\alpha_1(p)$ be the first eigenvalue of \eqref{L_up} and let $\varphi_1$ be the first eigenfunction of $-\frac{d^2}{d t^2}$ in $I$ with Dirichlet boundary condition, namely $\varphi_1(t)=\cos\left(\frac{\pi}{2}t\right)$ for $t\in\bar I$, then we have\vspace{0.2cm}
		
		\begin{enumerate}
			\item \begin{equation}\label{convergence_from_above}
				\Vert u_{p} \Vert_{\infty}^{p-1}=\frac{\pi^2}{4}+\frac{\pi^2}{4}\tilde c (p-1)+o(p-1)\qquad\text{as $p\to1$},
			\end{equation} 
			where
			\begin{equation*}
				\tilde{c}=\frac{\int_{-1}^{1}\varphi_1^2(t)|\log\left(\varphi_1(t)\right)|\,dt}{\int_{-1}^{1}\varphi_1^2(t)\,dt}>0;
			\end{equation*}
			\item\begin{equation*}
				\lim_{p\to1}\frac{u_p}{\|u_p\|}_{\infty}=\varphi_1\qquad \text{in $C^1(\bar I)$};
			\end{equation*}
			\item \begin{equation*}
				pu_p^{p-1}\to\frac{\pi^2}{4}\qquad\text{in $C^1_{loc}(I)$\qquad\text{as $p\to1$};}
			\end{equation*}
			\item \begin{equation*}
				\lim_{p\to1}\alpha_1(p)=0.
			\end{equation*}
		\end{enumerate}
	\end{teorema}
	The results concerning $p\to+\infty$ are inspired by \cite{Grossi1, GGPS}, where the behavior of solutions to the Lane-Emden problem in annuli of $\R^N$, $N\geq2$, is analyzed. A Liouville limit problem has been detected also in the asymptotic analysis, as $p\to+\infty$ of positive, finite energy solutions of the Lane-Emden problem in planar domains (see \cite{RenWeiTAMS1994, RenWeiPAMS1996, AdiGrossi, DIPpositive, DGIPsqrte, Thizy}).\\
	For what concerns the case $p\to1$, our analysis improves the one in \cite{Grossi2} in the one-dimensional problem in the spirit of what has been done in the planar unit ball in \cite{GladIan}.
	
	These theorems allow to prove our stability/instability results.

	\begin{teorema}[Stability for $p \to +\infty$]\label{Teo_stability_p_infty}
		Let $\omega \subset \mathbb{R}^{N-1}$ be a smooth bounded domain, let $\lambda_1(\omega)$ be the first non-trivial eigenvalue of the Laplacian in $\omega$ (i.e. in $(N-1)$ coordinates) with Neumann boundary condition, let $u_{p,L}$ be the positive one-dimensional solution to \eqref{problem_Omega} in $\Omega_{\varphi_L}$ defined in \eqref{upLN} and let $\alpha_1(p)$ be the first eigenvalue of $\eqref{L_up}$ with Dirichlet boundary conditions in I. Then, {for any $\gamma>\frac12$, if $p$ is sufficiently large and $L>\frac{\gamma\, p}{\sqrt{\lambda_1(\omega)}}$, the pair $(\Omega_{\varphi_L}, u_{p,L})$ is a stable energy-stationary pair (see Definition $\ref{def:critical_energy_pair_stable}$).}
	\end{teorema}
	\begin{teorema}[Stability/Instability for $p \to 1$]\label{Teo_stability_p_1}
		Let $\omega \subset \mathbb{R}^{N-1}$ be a smooth bounded domain, let $\lambda_1(\omega)$ be the first non-trivial eigenvalue of the Laplacian in $\omega$ (i.e. in $(N-1)$ coordinates) with Neumann boundary condition and let $u_{p,L}$ be the positive one-dimensional solution to \eqref{problem_Omega} in $\Omega_{\varphi_L}$ defined in \eqref{upLN}. Then:\\
		-- if $L<\sqrt{\frac{\pi^2}{4\lambda_1(\omega)}}$, for $p$ sufficiently close to $1$, the pair $({\Omega_{\varphi_L}}, u_{p,L})$ is an unstable energy-stationary pair  (see Definition $\ref{def:critical_energy_pair_unstable}$);\\
		-- if $L>\sqrt{\frac{\pi^2}{4\lambda_1(\omega)}}$, for $p$ sufficiently close to $1$, the pair $({\Omega_{\varphi_L}}, u_{p,L})$ is a stable energy-stationary pair (see Definition $\ref{def:critical_energy_pair_stable}$).
	\end{teorema}
	\begin{osservazione}
		When stability holds we have that the cylinder with the one-dimensional solution there definite locally minimizes the energy functional with the volume constrained. This does not hold when $(\Omega_{\varphi_L},u_{\varphi_L})$ is unstable, hence there are domains close to the bounded cylinder and with the same volume which are better candidates to optimize the energy.
	\end{osservazione}
	The proofs of the above theorems are based on an important general characterization of the stability/instability proved in \cite{AfonsoPacella} (see Theorem $\ref{teorema}$).\\
	As compared with the results of \cite{AfonsoPacella} we substantially improve them in the case of the Lane-Emden nonlinearity. This is because we precisely estimate the behaviour of the quantity $p\|u_p\|_{\infty}^{p-1}=pu_p^{p-1}(0)$ for $p$ large or for $p$ close to 1. Note that for $p$ close to 1 we have an explicit threshold for the instability for which there were no previous results in the nonlinear case. It is worth mentioning that another important ingredient in the proof is to have obtained in Theorem $\ref{Teo_p_to_inf}$ and Theorem $\ref{Teo_p_to_1}$ the asymptotic behaviour of the first eigenvalue $\alpha_1(p)$ of problem \eqref{L_up}. \\
	Our results have been given in the context of positive solutions of Lane-Emden problems. It would be interesting to study similar questions for sign-changing solutions and for other important semilinear problems as, for example, the Henon type ones.\\
	Finally, we recall that the study of the critical pairs of the energy functional $T$ is equivalent to the study of domains for which the overdetermined problem (with an extra homogeneous Neumann condition on the relative boundary) admits a solution, see \cite[Proposition 2.6]{AfonsoPacella}.\\
	Thus our instability result suggests that, for $p$ close to $1$, non trivial domains admitting solutions of the relative overdetermined problem should bifurcate from the bounded cylinder. Results in this direction have been proved in the case of the torsion problem in \cite{PacRuSic} and for an eigenvalue problem in \cite{last}.
	\\\\
	The paper is organized as follows. In Section $\ref{Section:Nondeg_asymt_behav}$ we study the asymptotic behaviour of $u_p$ as $p \to \infty$ and $p \to 1$ and prove Theorem $\ref{Teo_p_to_inf}$ and Theorem $\ref{Teo_p_to_1}$. In Section $\ref{section:stability}$ we first set clearly the question of the stability/instability, in the framework of domains which are hypographs, clarifying better the definition given in \cite{AfonsoPacella}. Then we prove Theorem $\ref{Teo_stability_p_infty}$ and Theorem $\ref{Teo_stability_p_1}$ by showing how the values of $p$ allow to satisfy the characterization given in \cite{AfonsoPacella} (see Theorem \ref{teorema}).

	\section{Nondegeneracy and asymptotic behaviour of solutions to \eqref{up_problem}}\label{Section:Nondeg_asymt_behav}
	We first show the nondegeneracy of the unique positive solution to \eqref{up_problem}, then in Subsection \ref{subsection:p_to_inf} and in Subsection \ref{subsection:p_to_1} we prove Theorem \ref{Teo_p_to_inf} and Theorem \ref{Teo_p_to_1} respectively.\\
	The following result is probably known, but for reader's convenience we provide the proof.
	\begin{proposizione}\label{prop:nondeg}
		For any $p>1$ and any $L>0$ the solution $u_{p,L}$ to \eqref{up_problem} is non-degenerate.    
	\end{proposizione}
	\begin{proof}
		We need to show that the eigenvalue problem
		\begin{equation*}
			\begin{cases}
				-z_p''-pu^{p-1}_{p,L}z_p=\alpha z_p \quad \mbox{ in } (0,L)\\
				z_p'(0)=z_p(L)=0
			\end{cases}
		\end{equation*}
		doesn't admit zero as an eigenvalue.
		Let us consider the even extension of both $u_{p,L}$ and $z_p$ to $(-L,L)$ (we still denote them by $u_{p,L}$ and $z_p$). If by contradiction there exists a solution $z_p\not\equiv 0$ of the problem 
		\begin{equation*}
			\begin{cases}
				-z_p''=pu^{p-1}_{p,L}z_p \quad \mbox{ in } (-L,L)\\
				z_p(\pm L)=0,
			\end{cases}
		\end{equation*}
		then, being $u_p$ the least energy solution (and hence with Morse index one), $z_p$ must be a second eigenfunction and have only two nodal regions. This is impossible since $z_p$ is even. Hence $z_p \equiv 0$.
	\end{proof}

	\subsection{Asymptotic behaviour as $p\to+\infty$}\label{subsection:p_to_inf}
	
	In this section, we prove Theorem $\ref{Teo_p_to_inf}$, in particular $(1)$ follows from Lemma \ref{lemma1}, $(2)$ from Theorem \ref{teogreen}, $(3)$ from Proposition \ref{prop:Lp_estimate}, $(4)$ from Proposition \ref{prop:convergencetoW} and $(5)$ from Theorem \ref{thm:alpha1p}.
	\begin{lemma}\label{lemma1}
		Let $u_p$ be the positive solution to \eqref{problema_in_I}, then we have that
		\begin{equation*}
			\Vert u_p \Vert_{\infty}^{p-1} \geq \nu_1(I)=\frac{\pi^2
			}{4} \quad \mbox{ for all } p>1\qquad\text{and}\qquad
			\Vert u_p \Vert_ \infty \rightarrow 1 \quad \mbox{as } p \to \infty,
		\end{equation*}
		where $\nu _1(I)$ is the first eigenvalue of $-\frac{d^2}{dx^2}$ in $I$ with Dirichlet boundary condition.
		\begin{proof}
			Integrating equation \eqref{problema_in_I} against $\varphi_1$, the first eigenfunction of $-\frac{d^2}{dx^2}$ in $I$ with Dirichlet boundary condition, we have
			\[
			\int_{-1}^{1}u^p_p(t)\varphi_1(t) \,dt=\int_{-1}^{1} -u_p''(t)\varphi_1(t)\,dt
			=\int_{-1}^{1} -u_p (t)\varphi''_1(t)\,dt 
			=\frac{\pi^2
			}{4}\int_{-1}^{1} u_p(t) \varphi_1 (t)\,dt
			\]
			that is
			\begin{equation*}
				\int_{-1}^{1} u_p(t) \varphi_1(t)(u^{p-1}_p(t)-\frac{\pi^2
				}{4})\,dt=0.
			\end{equation*}
			Thus, being both $u_p$ and $\varphi_1$ positive, then $u_p^{p-1}-\frac{\pi^2
			}{4}$ changes sign, clearly this implies that $\Vert u_p \Vert_{\infty}^{p-1} \geq \frac{\pi^2}{4}$ for all $p$, which in turn allows to deduce that 
			\begin{equation}\label{limsup_geq_1}
				\liminf_{p \to \infty}\Vert u_p \Vert_{\infty}\geq 1. \end{equation} 
			Next we look for a bound from above of $\Vert u_p \Vert_{\infty}$.\\
			Being $u_p$ concave, the graph of $u_p$ is above the graph of the piecewise linear function $g_p(t):=\Vert u_p\Vert_\infty(1-|t|)$ for any $t\in I$. In particular, if we fix $\alpha\in(0,1)$, then 
			\[
			u_p^p(t)\geq \Vert u_p\Vert ^p_\infty(1-\alpha)^p \qquad\text{for any $t\in[0,\alpha].$}
			\]
			Thus
			\[
			u'_p(t)=-\int_0^t u_p^p(\tau)d\tau\leq -t\Vert u_p\Vert^p_\infty(1-\alpha)^p\qquad\text{for ant $t\in[0,\alpha]$}
			\]
			and, as a consequence,
			\begin{eqnarray*}
				\Vert u_p\Vert_\infty(1-\alpha)&=&g_p(\alpha)\leq u_p(\alpha)= u_p(0)+\int_0^\alpha u'_p(\tau)d\tau\\
				&\leq& \Vert u_p\Vert_\infty(1-\Vert u_p\Vert_\infty^{p-1}(1-\alpha)^p\frac{\alpha^2}{2}).
			\end{eqnarray*}
			Hence
			\[
			\Vert u_p \Vert_{\infty} \leq  \frac{2^{\frac{1}{p-1}}}{\alpha^{\frac{1}{p-1}}(1-\alpha)^{\frac{p}{p-1}}} \xrightarrow[p\to \infty]{} \frac{1}{1-\alpha}.
			\]
			By the arbitrariness of $\alpha\in(0,1)$ we can deduce that
			\[ 
			\limsup_{p \to \infty} \Vert u_p \Vert_{\infty}\leq 1. \]
			Combining the latter estimate with \eqref{limsup_geq_1} we finally get
			\begin{equation*}
				\lim_{p \to \infty} \Vert u_p \Vert_{\infty} = 1.
			\end{equation*}
		\end{proof}
	\end{lemma}
	\begin{lemma}\label{lemma:Ip}
		Let $p>1$ and let us denote by
		\begin{equation}\label{Ip}
			I_p := \inf_{u \in H^1_{0}(I)} \frac{\int_{I}| u'(t)| ^2\,dt}{\biggl(\int_{I}u^{p+1}(t)\,dt\biggr)^{\frac{2}{p+1}}}.
		\end{equation}
		Then 
		\begin{equation*}
			\limsup_{p \to \infty} I_p \leq 2.
		\end{equation*}
		\begin{proof}
			We denote by $\omega(t)=1-|t|$, $t \in I$; then
			\[
			I_p \leq \frac{\int_{I}1\,dt}{\biggl(\int_{I}(1-|t|)^{p+1}\,dt\biggr)^{\frac{2}{p+1}}}= \frac{2}{2^{\frac{2}{p+1}}\biggl(\int_{0}^{1}(1-t)^{p+1}\,dt\biggr)^{\frac{2}{p+1}}}=2^{\frac{p-1}{p+1}}(p+2)^{\frac{2}{p+1}} 
			\]
			for all $p>1$, so the claim follows easily.
		\end{proof}
	\end{lemma}
	\begin{lemma}\label{lem:bounds}
		Let $u_p$ be the positive solution to \eqref{problema_in_I}, then there exist $c,\,C>0$ such that
		\begin{equation*}
			c\leq \int_{I} (u'_p(t))^2 \,dt=\int_{I} u^{p+1}_p(t)\,dt\leq C\qquad\text{as $p\to+\infty$.}
		\end{equation*}
		\begin{proof}
			Let $\tilde{u}_p\geq0$ be a minimizer of $I_p$, then $\tilde{u}_p$ solves
			\begin{equation*}
				-\tilde{u}''_p= \alpha _p \tilde{u}^p_p\quad\text{in $I$}, \qquad \mbox{ where } \alpha_p = \frac{\int_I (\tilde{u}'_p(t))^2\,dt}{\int_I \tilde{u}^{p+1}_p(t)\,dt},
			\end{equation*}
			then $u_p=\alpha_p^{\frac{1}{p-1} }\tilde{u}_p$ and
			\begin{equation*}
				\int_I (u_p'(t))^2\,dt =\int_{I} u^{p+1}_p(t)\,dt=\alpha ^{\frac{p+1}{p-1}} \int_{I} \tilde{u}^{p+1}_p(t)\,dt= I_p^{\frac{p+1}{p-1}},
			\end{equation*}
			where $I_p$ is defined in \eqref{Ip}.
			The boundedness from above then follows from Lemma \ref{lemma:Ip}.
			On the other hand, by Lemma \ref{lemma1} and the Rellich-Kondrachov theorem there exists $c>0$ such that
			\[
			1\leq\liminf_{p\to+\infty}\Vert u_p\Vert_\infty\leq c^{-1}\int_{I}(u'_p(t))^2dt= c^{-1}\int_{I} u^{p+1}_p(t)\,dt
			\]
		\end{proof}
	\end{lemma}
	
	\begin{lemma}
		Let $u_p$ be the positive solution to \eqref{problema_in_I}, then there exists $p_0>1$ such that
		\begin{equation*}
			\| u'_p \|_{\infty} \leq C \quad \mbox{ for all } p\geq p_0.
		\end{equation*}
		\begin{proof}
			The thesis can be easily obtained using the symmetry and the positivity of $u_p$ and applying  H\"older inequality and Lemma \ref{lem:bounds}
			\[
			\vert u_p'(t) \vert \leq \left|\int_{0}^{t} u_p^p(\tau)\,d\tau\right| \leq \left(\int_{0}^{1} u_p^{p+1}(\tau)\,d\tau\right)^{\frac{p}{p+1}} \leq C.
			\]
		\end{proof}
		\label{lemma2}
	\end{lemma}
	
	\begin{teorema}\label{teogreen}
		Let $u_p$ be the positive solution to \eqref{problema_in_I}, then 
		\begin{equation*}
			u_p \rightarrow 2G(\cdot,0)=1-|\cdot|\qquad \text{in $C^0([-1,1])\;$ and in $\;C^1_{loc}([-1,1]\setminus\lbrace 0 \rbrace)$,}
		\end{equation*}
		where $G$ is the Green's function of $-\frac{d^2}{d t^2}$ with Dirichlet boundary conditions defined in
		\eqref{Green}. 
		\begin{proof}
			By virtue of Lemma $\ref{lemma2}$ and Ascoli-Arzelà theorem there exists $\bar{u} \in C^0(\bar{I})$ such that, up to a subsequence, $u_p \rightarrow \bar{u}$, in $C^0(\bar{I})$ and from Lemma $\ref{lemma1}$ we have $\bar{u} \not\equiv 0$.\\
			Since $\Vert u_p\Vert_{\infty} \rightarrow \Vert \bar{u} \Vert_{\infty}$, then by Lemma \ref{lemma1} 
			\begin{equation*}
				\Vert u_p\Vert_{\infty}=u_p(0) \rightarrow 1=\bar{u}(0)=\Vert \bar{u} \Vert_{\infty}.
			\end{equation*}
			Next we claim that 
			\begin{equation}
				\bar{u}(t)<1\quad \mbox{for all $t \neq 0$}.
				\label{claim:ubar}
			\end{equation}
			
			We assume by contradiction that there exists $ \tilde{t} \neq 0$ such that $\bar{u}(\tilde{t})=1$. Without loss of generality by the evenness of $\bar{u}$, we can suppose $\tilde{t}>0$. Since $\bar{u}$ inherits the monotonicity of $u_p$ in $[0,\tilde{t}]$, then $\bar{u}(t)\equiv 1$ in $[0, \tilde{t}]$.\\ Furthermore, there exists $p_0$ such that for any $p\geq p_0$
			\begin{equation}\label{less1}
				u_p(\tfrac{\tilde{t}}{2}) <1,
			\end{equation}
			indeed, if this is not the case, there exists $p_n\to+\infty$ such that
			\(
			u_{p_n}(\tfrac{\tilde{t}}{2})\geq1.
			\)
			Thus, being $u_{p_n}$ decreasing in $[0,1]$ and by Lemma \ref{lemma1},
			\begin{eqnarray*}
				1\leq u_{p_n}(\tfrac{\tilde{t}}{2})&=&\|u_{p_n}\|_\infty+\int_{0}^{\tfrac{\tilde{t}}{2}}u'_{p_n}(s)\,ds\\
				&=&\|u_{p_n}\|_\infty-\int_0^{\tfrac{\tilde{t}}{2}}\int_0^s u_{p_n}^{p_n}(\tau)\,d\tau\, ds\\
				&\leq& \|u_{p_n}\|_\infty-\frac{\tilde{t}^2}{8}\to1-\frac{\tilde{t}^2}{8}<1\quad\text{as $n\to+\infty$},
			\end{eqnarray*}
			which is impossible.\\
			Then for $p\geq p_0$, \eqref{less1} and the monotonicity of $u_p$ in $[0,1]$ imply that 
			\begin{equation*}
				\vert u''_p(t)\vert= u_p^p(t)<1 \quad \mbox{ in } \left[\tfrac{\tilde{t}}{2}, \tilde{t}\right],
			\end{equation*}
			so, up to a subsequence, $u_p$ converges to $\bar{u}$ in $C^1\left(\left[\frac{\tilde{t}}{2}, \tilde{t}\right]\right)$ and passing to the limit in the equation in the weak form we have
			\begin{equation*}
				-\bar{u}''=1 \quad \mbox{ in } \left[\tfrac{\tilde{t}}{2}, \tilde{t}\right],
			\end{equation*}
			which is a contradiction against $\bar{u} \equiv 1$ in $\left[\frac{\tilde{t}}{2}, \tilde{t}\right]$.
			This concludes the proof of \eqref{claim:ubar}.\\
			In turn, by \eqref{claim:ubar}, for every $\delta>0$ there exists $p_\delta>1$ such that, for every $p \geq p_\delta$,
			\[\vert u''_p(t) \vert =u_p^p(t)<1 \quad \mbox{ for any } t\in[ -1, -\delta] \cup [\delta, 1],
			\]
			then, up to a subsequence, $u_p \rightarrow \bar{u}$ in $C^1 ([ -1, -\delta] \cup [\delta, 1])$
			and passing to the limit in the weak equation solved by $u_p$ we obtain that $\bar{u}$ solves
			\begin{equation*}
				\begin{cases}
					-\bar{u}''=0 \qquad \mbox{ in } [ -1, -\delta] \cup [\delta, 1 ] \\
					\bar{u}(\pm 1)=0.
				\end{cases}
			\end{equation*}
			By the arbitrariness of $\delta$ and being $\bar{u}(0)=1$, we can conclude that $\bar{u} = 2G(\cdot, 0)$ in $I$.
		\end{proof}
	\end{teorema}
	
	Next we estimate the p-th power of the $L^{\infty}$ norm of the solution.
	\begin{proposizione}\label{prop:Lp_estimate}
		Let $u_p$ be the positive solution to \eqref{problema_in_I}, then 
		\begin{equation*}
			\Vert u_p \Vert^{p+1}_{\infty}=\frac{p}{2}(1+o_p(1))
		\end{equation*}
	\end{proposizione}
	\begin{proof}
		Multiplying the equation solved by $u_p$, integrating it against $u'_p$ in $(0,1)$, setting $\bar u(t):=2G(t,0)$ for $t\in[0,1]$ and applying Theorem \ref{teogreen} we get
		\[
		\frac{\Vert u_p \Vert^{p+1}_{\infty}}{p+1} = \frac{(u'_p(1))^2}{2} \xrightarrow[p \to \infty]{} \frac{(\bar u'(1))^2}{2}=\frac{1}{2}.
		\]
	\end{proof}
	Now we prove the convergence of a suitable rescaling of $u_p$ around the origin to a solution of $-W''=e^W$ in $\R$.
	\begin{proposizione}\label{prop:convergencetoW}
		Let $\mu_p=\left(p\Vert u_p\Vert_\infty^{p-1}\right)^{-\frac{1}{2}}$ and let $\tilde u_p$ and $W$ be the functions defined in \eqref{utildep} and in \eqref{W} respectively,
		then 
		\begin{equation*}
			\tilde{u}_p\rightarrow W \quad \mbox{ in } C^1_{loc}(\mathbb{R}).
		\end{equation*}
		\begin{proof}
			By direct computations we have that $\tilde{u}_p$ solves 
			\begin{equation}\label{utilde_eq}
				\begin{cases}
					-\tilde{u}''_p=\left(1+\frac{\tilde{u}_p}{p}\right)^p \quad \mbox{ in } \left[ -\frac{1}{\mu_p}, \frac{1}{\mu_p}\right]\\
					\tilde{u}'_p(0)=0\\
					\tilde{u}_p(0)=0.
				\end{cases}
			\end{equation}
			Being $-p \leq \tilde{u}_p\leq 0$, we deduce that $0\leq |\tilde{u}''_p|\leq 1$, thus in any compact subset of $\mathbb{R}$ $\tilde{u}'_p$ is uniformly bounded and then Ascoli-Arzelà theorem implies that up to a subsequence $\tilde{u}_p \rightarrow W$ in $C^1_{loc}(\mathbb{R})$, where $W$ solves \eqref{W2}.
			It is easy to see that the solution of $\eqref{W2}$ takes the form 
			\begin{equation*}
				W(s)=\log\frac{4e^{\sqrt{2}s}}{(1+e^{\sqrt{2}s})^2}
			\end{equation*}
			and this concludes the proof.
		\end{proof}
	\end{proposizione}
	\begin{lemma}\label{lemma:decay}
		Let $\tilde u_p$ be the rescaled function defined in \eqref{utildep} and let $W$ be the function in \eqref{W}. Then there exist $C_1,\,C_2>0$ such that, for $p$ large enough,
		\begin{equation*}
			\tilde{u}_p \leq C_1 W+C_2 \quad \mbox{ in } \left [-\frac{1}{\mu _p}, \frac{1}{\mu _p} \right].
		\end{equation*}
		\begin{proof} It can be easily seen that $\int_{0}^{+\infty} e^{W(s)}\,ds = \sqrt{2}$, therefore there exists $R>0$ such that
			\begin{equation}
				\int_{0}^{R}e^{W(s)} \,ds \geq \frac{\sqrt{2}}{2}.
			\end{equation} 
			We know from Proposition $\ref{prop:convergencetoW}$ that
			\begin{equation}\label{utilde_to_W}
				\tilde{u}_p \rightarrow W \mbox{ in } C^1_{loc}(\mathbb{R}), 
			\end{equation}
			then for $p$ large enough 
			\begin{equation*}
				\tilde{u}_p \leq W+1 \quad \mbox{ in }[-R,R]. 
			\end{equation*}

			Now we have to prove the estimate in $\left[-\frac{1}{\mu _p},-R \right) \cup \left(R, \frac{1}{\mu _p} \right]$.\\ 
			By symmetry it is obviously enough to show it for $s\in \left(R, \frac{1}{\mu _p} \right]$.\\
			By \eqref{utilde_eq} we have
			\begin{equation*}
				\tilde{u}'_p(s)=-\int_{0}^{s} \left( 1+ \frac{\tilde{u}_p(r)}{p} \right)^p\,dr \leq -\int_{0}^{R} \left( 1+ \frac{\tilde{u}_p(r)}{p} \right)^p\,dr,
			\end{equation*}
			so by \eqref{utilde_to_W}, as $p\to+\infty$
			\begin{equation}\label{stima_u'p}
				\tilde{u}'_p(s) \leq -\int_{0}^{R} e^{W(r)} \,dr {+o_p(1)}\leq -\frac{\sqrt{2}}{4}.
			\end{equation}
			In conclusion we get the thesis being
			\begin{equation*}
				\tilde{u}_p(s) =\tilde{u}_p(R) +\int_{R}^{s}\tilde{u}'_p(r)\,dr\leq -\frac{\sqrt{2}}{4}(s-R)\leq \frac{W(s)}{4}+\frac{\log2}{2}+\frac{\sqrt2}{4}R,
			\end{equation*}
			where the inequalities follow from $\tilde u_p\leq0$, \eqref{stima_u'p} and 
			\[
			W(s)=\log \frac{4e^{-\sqrt{2}s}}{(1+e^{-\sqrt{2}s})^2}=2\log 2 -\sqrt{2}s-2\log(1+e^{-\sqrt{2}s})\geq -\sqrt{2}s-2\log2.
			\]
		\end{proof}
	\end{lemma}

	\begin{proposizione}\label{prop:alpha1pmup2_to_beta}
		Let $\alpha_1(p)$ be the first eigenvalue of \eqref{L_up}, then, up to subsequences, we have
		\begin{equation}\label{def_beta_1}
			\alpha _1(p) \mu_p ^2 \xrightarrow[p \to \infty]{} \beta _1<0.
		\end{equation}
		\begin{proof}
			As a first step we show that for $p$ large enough
			\begin{equation}\label{bounds_alpha1p}
				-1< \alpha_1(p)\mu_p^2<0.
			\end{equation}
			
			Let us consider the following test function
			\[
			\phi (t)=
			\begin{cases}
				\tfrac{\sqrt3}{2} &\quad \mbox{ if } t \in [-\frac12, \frac{1}{2}]\\
				\sqrt3 (1-|t|) &\quad \mbox{ if } t \in I\setminus[-\frac12, \frac{1}{2}].
			\end{cases}
			\]
			{ Being $\Vert \phi \Vert_2=1$, the variational characterization of eigenvalues, direct computations, Lemma \ref{teogreen}, combined with the fact that $2G(t,0)=1-t\leq \frac12$ for any $t\in[\frac12,1]$, Lemma \ref{lemma1} and Lemma \ref{lem:bounds} imply that
				\begin{eqnarray}\label{alpha1p_alto}
					\alpha_1(p)&=& \inf _{v \in H^1_0(-1,1)}\frac{\int_{-1}^{1}|v'(t)|^2\,dt-\int_{-1}^{1}pu_p^{p-1}(t)v^2(t)\,dt}{\int_{-1}^{1}v^2(t)\,dt}\nonumber\\
					&\leq&3-\int_{-1}^{1}pu_p^{p-1}(t)\phi^2(t)\,dt\leq 3-\frac34\int_{-\frac12}^{\frac12}pu_p^{p-1}(t)\,dt\nonumber\\
					&\leq&3-\frac34\int_{-1}^{1}pu_p^{p-1}(t)\,dt+\frac32\int_{\frac12}^{1}pu_p^{p-1}(t)\,dt\nonumber\\
					&\leq&3-\frac34 \frac{p}{\Vert u_p\Vert_\infty^2}\int_{-1}^{1}u_p^{p+1}(t)\,dt+o_p(1)\nonumber\\
					&\leq&3-\frac34 \frac{p}{1+o_p(1)} c+o_p(1)
					<0 \quad \mbox{as } p \to +\infty
				\end{eqnarray}
			}
			
			On the other hand let $w_{1,p}>0$ the eigenfunction associated to $\alpha_1(p)$ satisfying $\Vert w_{1,p}\Vert_{L^2(I)}=1$, then 
			\begin{equation}\label{alpha1p_basso}
				\alpha_1(p)= \frac{\int_{-1}^{1}|w_{1,p}'(t)|^2\,dt-\int_{-1}^{1}pu_p^{p-1}(t)w_{1,p}^2(t)\,dt}{\int_{-1}^{1}w_{1,p}^2(t)\,dt}> -p\Vert u_p\Vert_\infty^{p-1}. 
			\end{equation}
			In conclusion \eqref{bounds_alpha1p} follows combining \eqref{alpha1p_alto} with \eqref{alpha1p_basso}.\\
			Hence, in order to conclude the proof it is enough to exclude that $\alpha_1(p)\mu_p^2\to0$.\\ 
			If we assume by contradiction that $\alpha _1(p) \mu_p ^2 \rightarrow 0$, then, setting $\lambda_p :=-\alpha_1(p) +1$, we have 
			\begin{equation}\label{lambdapmup2_to_0}
				\lambda_p + \alpha_1(p) >0 \qquad \mbox{ and } \qquad\lambda_p \mu^2_p \rightarrow 0.
			\end{equation}
			So, by maximum principle, we can deduce that any solution $k_p$ of 
			\begin{equation*}
				\begin{cases}
					k_p''=(\lambda_p - pu^{p-1})k_p \quad \mbox{ in $I$} \\
					k_p( \pm 1)= \vert u'_p(1) \vert
				\end{cases}
			\end{equation*}
			is positive and even in $I$. Furthermore, being $k''_p(0)=p\|u_p\|_\infty^{p-1}(\lambda_p\mu_p^2-1)k_p(0)<0$ for $p$ large enough and $k'_p(0)=0$:
			\begin{itemize}
				\item[\emph{i)}] either $\Vert k_p\Vert_\infty=k_p(0)$,
				\item[\emph{ii)}] or, if the above condition is not fulfilled, there exists $m_p\in(0,1)$ such that $k_p$ is decreasing in $(0,m_p)$ and $k_p$ is increasing in $(m_p,1)$.
			\end{itemize} 
			In the latter case $m_p>t_p$, where $t_p$ is the positive inflection point of $k_p$, namely $t_p\in(0,1)$ such that $\lambda_p=p u_p^{p-1}(t_p)$. Let us show that 
			\begin{equation}\label{mp_mup}
				\frac{m_p}{\mu_p}\to+\infty.
			\end{equation}
			If by contradiction \eqref{mp_mup} does not hold, then, up to a subsequence $\frac{t_p}{\mu_p}\to t_\infty\in [0,+\infty)$, then, by Proposition \ref{prop:convergencetoW},
			\begin{equation*}
				0 \leftarrow \lambda_p\mu^2_p= \frac{p u_p^{p-1}(t_p)}{p\Vert u_p\Vert_\infty^{p-1}}= \left( 1+\frac{\tilde{u}_p(\frac{t_p}{\mu_p})}{p}\right)^{p-1}\rightarrow \,e^{W({t}_{\infty})}>0.
			\end{equation*}
			From the contradiction above we get \eqref{mp_mup}.\\
			Now we set 
			\begin{equation*}
				y_p= 
				\left\{\begin{array}{ll}
					1 &\mbox{ if } \Vert k_p \Vert_{\infty}=k_p(0)\\
					m_p &\mbox{ otherwise }
				\end{array}\right. 
			\end{equation*}
			and we notice that, by \eqref{mp_mup}, in both cases $\tilde{y}_p:=\frac{y_p}{\mu_p} \rightarrow +\infty$ and  $k_p(0)= \Vert k_p \Vert_{L^{\infty}[-y_p, y_p]}$.\\
			The rescaled function
			\begin{equation*}
				\tilde{k}_p(s)=\frac{k_p(\mu_p s)}{k_p(0)}, \qquad y \in \left[ -\tilde{y}_p,\tilde{y}_p \right],
			\end{equation*}
			solves
			\begin{equation*}
				-\tilde{k}''_p (s) =\left( \mu_p^2 \lambda_p - \left(1+\frac{\tilde{u}_p(s)}{p}\right)^{p-1} \right) \tilde{k}_p(s),\qquad \tilde k_p(0)=1=\Vert \tilde k_p\Vert_\infty.
			\end{equation*}
			By \eqref{lambdapmup2_to_0} and observing that $\tilde u_p(s)\in [-p,0]$ for any $s\in[-\tilde{y}_p, \tilde{y}_p]$, in any compact subset $D$ of $[-\tilde{y}_p, \tilde{y}_p]$
			\[
			\vert \tilde{k}'_p(s) \vert \leq \int_{0}^{|s|}\left| \mu_p^2 \lambda_p - \left(1+\frac{\tilde{u}_p(r)}{p}\right)^{p-1} \right| \tilde{k}_p(r)\,dr  \leq \int_{0}^{|s|} (|o_p(1)| + 1) \,dr \leq C_D
			\]
			and 
			\[
			\vert \tilde{k}''_p(s) \vert \leq \left| \mu_p^2 \lambda_p - \left(1+\frac{\tilde{u}_p(s)}{p}\right)^{p-1} \right| \tilde{k}_p(s) \leq(|o_p(1)| + 1) \leq 2.
			\]
			Thus by Ascoli-Arzelà theorem, up to subsequences, $\tilde{k}_p \to\tilde{K}$ in $C^1_{loc}(\R)$, where $\tilde K$ solves
			\begin{equation*}
				\begin{cases}
					-\tilde{K}''=e^W \tilde{K} \quad \mbox{ in } \mathbb{R}\\
					\tilde{K}(0)=1, \quad\tilde{K}'(0)=0
				\end{cases}
			\end{equation*}
			and takes the form (see \cite[Lemma 4.2]{Grossi1})
			\begin{equation*}
				\tilde{K}(s)=1+\frac{s}{\sqrt{2}}\left(\frac{1-e^{\sqrt{2}s}}{1+e^{\sqrt{2}s}}\right).
			\end{equation*}
			It is immediate to see that $\tilde{K}(s)\to-\infty$ as $s\to+\infty$ and
			\[
			\tilde{K}'(s)= \frac{1-e^{\sqrt{2}s}}{\sqrt{2}(1+e^{\sqrt{2}s})}-\frac{2se^{\sqrt{2}s}}{(1+e^{\sqrt{2}s})^2}<0\qquad\text{for any $s>0$},
			\]
			then there exists a unique $\bar{R}>0$ such that $\tilde{K}(\bar{R})=0$ and $\tilde{K}(s)<0$ for any $s>\bar R$.
			Finally we reach the desired contradiction being
			\begin{equation*}
				0<\tilde{k}_p(2\bar{R)} \rightarrow \tilde{K}(2\bar{R})<0.
			\end{equation*}
		\end{proof}
	\end{proposizione}
	\begin{teorema}\label{thm:alpha1p}
		Let $\alpha_1(p)$ be the first eigenvalue of \eqref{L_up}, then 
		\[
		\alpha_1(p)={-\frac{1}{4}}p^2 (1+o_p(1)) \qquad\text{as} \quad p\to+\infty.
		\]
	\end{teorema}
	\begin{proof}
		Let $w_{1,p}$ be the eigenfunction associated to $\alpha_1(p)$ such that $\Vert w_{1,p} \Vert_{\infty}=1$ and $w_{1,p}>0$. It solves the following ODE problem
		\begin{equation*}
			\begin{cases}
				-w''_{1,p}=pu^{p-1}_pw_{1,p}+\alpha_1(p)w_{1,p} \quad &\mbox{ in } I\\
				w_{1,p}(\pm 1)=0.
			\end{cases}
		\end{equation*}
		The rescaled function $\tilde{w}_{1,p}(y)=w_{1,p}(\mu _py)$ solves
		\begin{equation}\label{eq:wtilde1p}
			-\tilde{w}''_{1,p}=\left( 1+\frac{\tilde{u}_p}{p}\right)^{p-1}\!\!\!\!\tilde{w}_{1,p}+\mu^2_p\alpha_1(p)\tilde{w}_{1,p}\qquad \text{in $I/\mu_p$,}
		\end{equation}
		where $\tilde u_p$ is defined in \eqref{utildep}, then, observing that $\tilde u_p(s)\in [-p,0]$ for any $s\in I/\mu_p$ and applying Proposition \ref{prop:alpha1pmup2_to_beta}, we get, for any $s\in I/\mu_p$:
		\begin{align*}
			|\tilde{w}'_{1,p}(s)| &\leq \int_{0}^{|s|}\left|\left( 1+\frac{\tilde{u}_p(r)}{p}\right)^{p-1}\!\!+\mu^2_p\alpha_1(p)\right|\tilde{w}_{1,p}(r)\,dr\\
			&\leq \int_{0}^{|s|} (1+\mu^2_p|\alpha_1(p)|)dr\leq C|s|\qquad\text{as $p\to+\infty$}
		\end{align*}
		and 
		\[
		|\tilde{w}''_{1,p}(s)| \leq 1+\mu^2_p|\alpha_1(p)|\leq 2.
		\]
		Thus, by Ascoli-Arzelà theorem,  $\tilde{w}_{1,p}\rightarrow \psi_1$ in $C^1_{loc}(\R)$ and $\psi_1$ is the solution to
		\begin{equation*}
			\begin{cases}
				-\psi''_1=e^W\psi_1+\beta_1\psi_1 \quad \mbox{ in }\mathbb{R}\\
				\psi_1\geq 0,
			\end{cases}
		\end{equation*}
		{where $\beta_1$ is the constant introduced in \eqref{def_beta_1}.}\\
		Let us first show that $\psi_1 \not\equiv 0$.\\
		Let $\xi _p$ be such that $w_{1,p}(\xi_p)=\Vert w_{1,p} \Vert_{\infty}=1$ and set $\tilde{\xi}_p=\frac{\xi_p}{\mu_p}$, then we have
		\begin{equation*}
			\tilde{w}_{1,p}(\tilde{\xi}_p)=1, \quad \tilde{w}'_{1,p}(\tilde{\xi}_p)=0, \quad \tilde{w}''_{1,p}(\tilde{\xi}_p)\leq 0.
		\end{equation*}
		If $|\tilde{\xi}_p|\rightarrow +\infty$ then, by \eqref{eq:wtilde1p}, Lemma \ref{lemma:decay}, Proposition \ref{prop:alpha1pmup2_to_beta} and \eqref{W}, for $p$ large enough,
		\begin{equation*}
			0\leq -\tilde{w}''_{1,p}(\tilde{\xi}_p)\leq e^{\frac{p-1}{p}(C_1W(\tilde{\xi}_p)+C_2)}+\beta_1+o_p(1)=\beta_1 + o_p(1)<0,
		\end{equation*}
		that is impossible. So we deduce that $|\tilde{\xi}_p|\leq C$ and, up to subsequences, $\tilde{\xi}_p\rightarrow\tilde{\xi}_{\infty}\in\R$ and $1=\tilde{w}_{1,p}(\tilde{\xi}_p)\rightarrow \psi_1(\tilde{\xi}_{\infty})$, which implies that $\psi_1 \neq 0$.\\
		Finally we compute $\beta_1$ setting $t=\sqrt{2}\log s$ and $Z(s)=\psi_1(t)$. Then $Z\neq 0$ and solves
		\begin{equation*}
			\begin{cases}
				-Z''-\frac{Z'}{s}=\frac{8}{(1+s^2)^2}Z+\frac{2\beta_1}{s^2}Z \quad &\mbox{ in }(0,+\infty)\\
				Z\geq 0 \quad &\mbox{ in } (0,+\infty)\\
				\Vert Z \Vert_{\infty} \leq 1,
			\end{cases}
		\end{equation*}
		whose solutions have the form (see \cite{GGPS})
		\begin{equation*}
			Z(s)=
			\begin{cases}
				c_1\frac{s^{\sqrt{-2\beta_1}}}{1+s^2}\left( s^2+\frac{\sqrt{-2\beta_1}+1}{\sqrt{-2\beta_1}-1}\right) +c_2\frac{s^{-\sqrt{-2\beta_1}}}{1+s^2}\left( s^2+\frac{\sqrt{-2\beta_1}-1}{\sqrt{-2\beta_1}+1}\right) &\mbox{ if } \beta_1 \neq -\frac{1}{2}\\
				c_1\frac{s}{1+s^2}+c_2\frac{s^4+4s^2\log s-1}{s(1+s^2)} &\mbox{ if } \beta_1=-\frac{1}{2},
			\end{cases}
		\end{equation*}
		therefore, being $Z$ bounded, we obtain $\beta_1=-\frac{1}{2}$.
	\end{proof}

	\subsection{Asymptotic behaviour as $p\to1$}\label{subsection:p_to_1}
	In this section, we prove Theorem $\ref{Teo_p_to_1}$, in particular $(1)$ follows from Proposition \ref{u_p(0)^{p-1}:close_to_1}, $(2)$ from Proposition \ref{u_p_norm_convergence}, $(3)$ from Corollary \ref{cor:Cloc_convergence} and $(4)$ from Proposition \ref{alpha_p:convergence}.
	\begin{proposizione}\label{u_p_norm_convergence}
		Let $u_p$ be the positive solution to \eqref{problema_in_I}, then 
		\begin{align}
			\label{norm_power_to_nu2}
			& \Vert u_p\Vert_{\infty}^{p-1}=u_p^{p-1}(0) \xrightarrow[p \to 1]{}\nu_1(I)=\frac{\pi^2}{4},\\
			& \bar u_p:=\frac{u_p}{\Vert u_p\Vert_{\infty}}\xrightarrow[p \to 1]{} \varphi_1\qquad \text{in $C^1(\bar I)$,}\label{baru_p}
		\end{align}
		where $\varphi_1(t)=\cos(\frac\pi2 t)$, $t\in I$, is the first eigenfunction and $\nu_1(I)$ is the first eigenvalue of $-\frac{d^2}{dt^2}$ in $I$ with Dirichlet boundary condition.
		\begin{proof}
			Let $\nu_p=\|u_p\|_\infty^{{p-1}}$, then $\bar{u}_p$, defined in \eqref{baru_p} is an even function solving
			\begin{equation}\label{baru_p equation}
				\begin{cases}
					-\bar u''_p=\nu_p\bar{u}_p^p \quad \mbox{ in } I\\
					\bar{u}_p(0)=1\\
					\bar{u}_p(\pm 1)=0.
				\end{cases}
			\end{equation}
			First of all we aim to show that $\nu_p$ is bounded.\\
			On the one hand
			\begin{equation}\label{stima_baru'_p}
				|\bar{u}'_p(t)|\leq \left|\int_{0}^{t}\nu_p\bar{u}_p^p(s)\,ds \right|\leq \nu_p|t|,\quad \text{for any $t\in I$,}
			\end{equation}
			then
			\begin{equation*}
				1=\bar{u}_p(0)\leq \int_{0}^{1}|\bar{u}'_p(t)|\,dt \leq \frac{1}{2}\nu_p
			\end{equation*}
			and so $\nu_p\geq 2$.\\
			On the other hand, if we have that up to subsequences $\nu_p \rightarrow +\infty$, then setting $w_p(s)=\bar{u}_p(\frac{s}{\sqrt{\nu_p}})$, $w_p$ solves
			\begin{equation*}
				\begin{cases}
					-w_p''=w_p^p\quad \mbox{ in } (-\sqrt{\nu_p}, \sqrt{\nu_p})\\
					w_p(0)=1\\
					w'_p(0)=0
				\end{cases}
			\end{equation*}
			and we have as before
			\begin{align*}
				w_p>0 \quad&\text{in $(-\sqrt{\nu_p}, \sqrt{\nu_p})$},\qquad w_p(0)=1,\\
				|w_p''| \leq w_p^p\leq 1,\quad \text{and} \quad |w'_p(s)|\leq &\left|\int_{0}^{s}|w''_p(r)|\,dr \right|\leq |s|,\quad \text{for any $s\in(-\sqrt{\nu_p}, \sqrt{\nu_p})$,}
			\end{align*}
			thus both $w_p$ and its first and second derivatives are bounded in compact sets. Thus, as a consequence of Ascoli-Arzelà theorem, we have that $w_p \rightarrow w$ in $C^1_{loc}(\mathbb{R})$ where $w$ solves
			\begin{equation}\label{problem_w}
				\begin{cases}
					-w''=w \quad \mbox{ in } \mathbb{R}\\
					w(0)=1\\
					w'(0)=0.
				\end{cases}
			\end{equation}
			Since the solution of \eqref{problem_w} is $w(s)=\cos(s)$, this is a contradiction against $w_p>0$.\\
			Thus we have that $\nu_p$ is bounded and, up to a subsequence, we denote by $\nu=\lim_{p \to 1} \nu_p$. Then, being $\|\bar u_p\|_\infty=1$, by \eqref{baru_p equation} and \eqref{stima_baru'_p}
			\[
			|\bar u''_p|\leq \nu_p \bar u_p^p\leq C\qquad\text{and}\qquad |\bar u_p'|\leq \nu_p\leq C\quad\:\text{in $I$ \quad\; for $p$ close to $1$.}
			\]

			This means that for Ascoli-Arzelà theorem, up to a subsequence, $\bar{u}_p\rightarrow \varphi$ in $C^1(\bar{I})$ and $\varphi$ solves
			\begin{equation*}
				\begin{cases}
					-\varphi''=\nu \varphi \quad &\mbox{ in } I\\
					\varphi>0 &\mbox{ in } I\\
					\varphi(\pm 1)=0.
				\end{cases}
			\end{equation*}
			This is an eigenvalue problem with Dirichlet boundary condition and $\nu$ must be the first eigenvalue, then $\nu=\nu_1(I)=\frac{\pi^2}{4}$ and $\varphi=\varphi_1$.
		\end{proof}
	\end{proposizione}
	\begin{corollario}\label{cor:Cloc_convergence}
		Let $u_p$ be the positive solution to \eqref{problema_in_I}, then 
		\begin{equation*}
			pu_p^{p-1}\to\frac{\pi^2}{4}\qquad\text{in $C_{loc}(I)$.}   
		\end{equation*}
	\end{corollario}
	\begin{proof}
		Let $\eta\in(0,1)$ and let $I_\eta=[-\eta,\eta]$, then 
		\begin{eqnarray*}
			\sup_{t\in I_\eta}\left|pu_p^{p-1}(t)-\frac{\pi^2}{4}\right|&=&\sup_{t\in I_\eta}\left|p\|u_p\|^{p-1}_\infty\bar u_p^{p-1}(t)-\frac{\pi^2}{4}\right|\\
			&\leq& \left|p\|u_p\|^{p-1}_\infty-\frac{\pi^2}{4}\right|\sup_{t\in I_\eta}|\bar u_p^{p-1}(t)|+\frac{\pi^2}{4}\sup_{t\in I_\eta}|\bar u_p^{p-1}(t)-1|\\
			&\leq&\left|p\|u_p\|^{p-1}_\infty-\frac{\pi^2}{4}\right|+\frac{\pi^2}{4}|\bar u_p^{p-1}(\eta)-1|\to0\quad\text{as $p\to1$},
		\end{eqnarray*}
		where $\bar u_p$ is defined in \eqref{baru_p} and we have used, in order to conclude, \eqref{norm_power_to_nu2} and that, by \eqref{baru_p}, $\bar u_p^{p-1}(\eta)=(\cos(\frac\pi2\eta))^{p-1}(1+o_p(1))$.\\
		By the arbitrariness of $\eta$ the thesis follows.
	\end{proof}
	\begin{proposizione}\label{u_p(0)^{p-1}:close_to_1}
		Let $u_p$ be the positive solution to \eqref{problema_in_I}, then 
		\begin{equation*}
			\|u_{p}\|^{p-1}_{\infty}=\frac{\pi^2}{4}+\tilde{c}\frac{\pi^2}{4}(p-1)+o(p-1),
		\end{equation*}
		where
		\begin{equation*}    \tilde{c}=\frac{\int_{-1}^{1}\cos^2\left(\frac{\pi}{2}t\right)|\log\left(\cos\left(\frac{\pi}{2}t\right)\right)|\,dt}{\int_{-1}^{1}\cos^2\left(\frac{\pi}{2}t\right)\,dt}>0.
		\end{equation*}
		\begin{proof}
			As noticed in the proof of Proposition $\ref{u_p_norm_convergence}$, $\bar u_{p}$ solves
			\begin{equation}\label{eq_v_p}
				\begin{cases}
					-\bar u_{p}''=\|u_{p}\|_\infty^{p-1}\bar u_{p}^{p-1} \quad \text{ in } I\\
					\bar u_{p}(0)=1\\
					\bar u_{p}(\pm 1)=0
				\end{cases}
			\end{equation}
			and $\bar u_p\to\varphi_1$ in $C^1(\bar I)$, where $\varphi_1(t)=\cos\left(\frac{\pi}{2}t\right)$ is the solution to 
			\begin{equation}\label{eq_varphi_1}
				\begin{cases}
					-\varphi''_1=\frac{\pi^2}{4}\varphi_1 \quad &\text{ in } I\\
					\varphi_1>0 \quad &\text{ in } I\\
					\varphi_1(\pm 1)=0.
				\end{cases}
			\end{equation}
			Now, multiplying $\eqref{eq_v_p}$ by $\varphi_1$, $\eqref{eq_varphi_1}$ by $\bar u_{p}$ and integrating both equations in $I$, we obtain
			\begin{equation}\label{identity_v_p_varphi}
				\|u_p\|_\infty^{p-1}\int_{-1}^{1}\bar u_{p}^{p}\varphi_1\,dt= \int_{-1}^{1}\bar u'_{p}\varphi_1'\,dt=\frac{\pi^2}{4}\int_{-1}^{1}\varphi_1\bar u_{p}\,dt.
			\end{equation}
			This implies that
			\begin{equation}\label{july}
				\left(\|u_p\|_\infty^{p-1}-\frac{\pi^2}{4} \right)\int_{-1}^{1}\bar u_{p}^{p}\varphi_1 \,dt=-\frac{\pi^2}{4}\int_{-1}^{1}\bar u_{p}\varphi_1\left(\bar u_{p}^{p-1}-1\right)\,dt.
			\end{equation}
			Next, using the identity $e^x-1=x \int_{0}^{1}e^{sx}\,ds$
			with $x=(p-1)\log \bar u_p$, we get
			\begin{equation*}
				\int_{-1}^{1}\bar u_p\varphi_1\left( \bar u_{p}^{p-1}-1\right)\,dt=\int_{-1}^{1}\varphi_1 \bar u_{p}(p-1)\log \bar u_{p}\int_{0}^{1}\bar u_{p}^{s(p -1)}\,ds \,dt,
			\end{equation*}
			which, combined with \eqref{july}, leads to
			\begin{equation*}
				\left(\|u_p\|_\infty^{p-1}-\frac{\pi^2}{4} \right)\int_{-1}^{1}\bar u_{p}^{p}\varphi_1 \,dt=\frac{\pi^2}{4}(p -1)\int_{-1}^{1}\varphi_1 \bar u_{p}\,|\log \bar u_p|\int_{0}^{1}\bar u_{p}^{s(p-1)}\,ds \,dt
			\end{equation*}
			which implies 
			\begin{equation*}
				\frac{\|u_p\|_\infty^{p-1}-\frac{\pi^2}{4} }{\frac{\pi^2}{4}(p-1)}=\frac{\int_{-1}^{1}\varphi_1 \bar u_{p}|\log \bar u_{p}|\int_{0}^{1}\bar u_{p}^{s(p -1)}\,ds \,dt}{\int_{-1}^{1}\bar u_{p}^{p}\varphi_1 \,dt},
			\end{equation*}
			what we want to show is that the right term of the last equality converges to an explicit constant $\tilde{c}>0$. The uniform convergence of $\bar u_{p}$ to $\varphi_1$ implies that
			\begin{equation*}
				\int_{-1}^{1}\bar u_p^{p}\varphi_1 \,dt \rightarrow\int_{-1}^{1}\varphi_1^2 \, dt > 0.
			\end{equation*}
			Moreover we know that $0 < \bar u_{p} \leq 1$ in $I$, therefore $\bar u_{p}\log \bar u_{p}$ is bounded and then 
			\begin{equation*}
				\|\varphi_1 \bar u_{p}|\log \bar u_{p}|\int_{0}^{1}\bar u_{p}^{s(p-1)}\,ds\|_{\infty} \leq \| \varphi_1 \bar u_{p}|\log \bar u_{p} \|_{\infty} \leq C.
			\end{equation*}
			Thus, using that $\bar u_{p} \to \varphi_1$ and the dominated convergence theorem 
			\begin{equation*}
				\int_{-1}^{1}\varphi_1 \bar u_{p}|\log \bar u_{p}|\int_{0}^{1}\bar u_{p}^{s(p-1)}\,ds \,dt \xrightarrow[p\to 1]{}\int_{-1}^{1}\varphi_1^2|\log \varphi_1| \,dt >0
			\end{equation*}
			then
			\begin{equation*}
				\frac{\|u_p\|_\infty^{p-1}-\frac{\pi^2}{4} }{\frac{\pi^2}{4}(p-1)}\longrightarrow\frac{\int_{-1}^{1}\varphi_1^2|\log \varphi_1| \,dt}{\int_{-1}^{1}\varphi_1^2 \,dt}=: \tilde{c}.
			\end{equation*}    
		\end{proof}
	\end{proposizione}

	\begin{proposizione} \label{alpha_p:convergence}
		Let $\alpha_1(p)$ be the first eigenvalue of \eqref{L_up}, then, up to subsequences, we have
		\begin{equation*}
			\alpha_1(p) \xrightarrow[p \to 1]{}0.
		\end{equation*}
		\begin{proof}
			For any $v\in H^1_0(I)$
			\[
			\frac{\int_{-1}^{1}|v'(t)|^2\,dt}{\int_{-1}^{1}|v(t)|^2\,dt}-p\|u_p\|_\infty^{p-1}\leq \frac{\int_{-1}^{1}|v'(t)|^2\,dt-p\int_{-1}^{1}u_p^{p-1}(t)v^2(t)\,dt}{\int_{-1}^{1}v^2(t)\,dt}\leq \frac{\int_{-1}^{1}|v'(t)|^2\,dt}{\int_{-1}^{1}|v(t)|^2\,dt},
			\]
			then, passing to the infimum on $v\in H^1_0(I)$, recalling the variational characterization of $\alpha_1(I)$, \eqref{norm_power_to_nu2} and that the first eigenvalue of $-\frac{d^2}{dt^2}$ in $I$ with Dirichlet boundary conditions is $\nu_1(I)=\frac{\pi^2}{4}$, we get
			\begin{equation}\label{bounds_alpha1p_pto1}
				o_p(1)\leq \alpha_1(p)\leq \nu_1(I)=\frac{\pi^2}{4}.
			\end{equation}
			Then, up to subsequences, we can define $\alpha:=\lim_{p\to1}\alpha_1(p)$.\\
			The eigenfunction $w_{1,p}$, associated to $\alpha_1(p)$ in $I$ and such that $\|w_{1,p}\|_\infty=1$ and $w_{1,p}>0$, solves
			\begin{equation*}
				\begin{cases}
					-w''_{1,p}-pu_p^{p-1}w_{1,p}=\alpha_1(p)w_{1,p} \quad &\mbox{ in }I\\
					w_{1,p}>0 &\mbox{ in } I\\
					w_{1,p}(\pm 1)=0.
				\end{cases}
			\end{equation*}
			By \eqref{bounds_alpha1p_pto1}, \eqref{norm_power_to_nu2} and being $w_{1,p}$ even, again, for any $t\in \bar I$ and for $p$ sufficiently close to $1$, we have 
			\begin{equation*}
				|w''_{1,p}(t)| \leq (\alpha_1(p)+pu_p^{p-1}(t))\| w_{1,p} \|_{\infty} \leq C,\quad\text{and}\quad |w'_{1,p}(t)|\leq \left|\int_{0}^{t}|w''_{1,p}(\tau)|\,d\tau\right|\leq C.
			\end{equation*}
			Thus, by Ascoli-Arzelà theorem, up to a subsequence, $w_{1,p}\to w_1$ in $C^1(\bar I)$ as $p\to1$ and, by Corollary \ref{cor:Cloc_convergence}, $w_1$ solves
			\begin{equation*}
				\begin{cases}
					-w''_1=\left( \frac{\pi^2}{4} + \alpha \right) w_1 \quad &\mbox{ in }(-1,1)\\
					w_1\geq 0 &\mbox{ in }(-1,1)\\
					w_1(\pm 1)=0.
				\end{cases}
			\end{equation*}
			Let $\xi_p \in (-1,1)$ be such that $w_{1,p}(\xi _p)=\Vert w_{1,p} \Vert_{\infty} =1$. Then $\xi_p$ is bounded so up to subsequences it converges to $\xi_1$ as $p\to1$. As a consequence, being
			\begin{align*}
				|1-w_1(\xi_1)|&=|w_{1,p}(\xi_p)-w_1(\xi_1)| \\
				&\leq |w_{1,p}(\xi_p)-w_{1,p}(\xi_1)|+|w_{1,p}(\xi_1)-w_1(\xi_1)|\\
				&\leq \|w_{1,p}'\|_{\infty}|\xi_p -\xi_1|+\|w_{1,p}-w_1 \|_{\infty}\\
				&\leq C|\xi_p -\xi_1|+\|w_{1,p}-w_1 \|_{\infty} \to 0 \quad \text{ as } p\to 1,
			\end{align*}
			we deduce that $w_1(\xi_1)=1$. Thus, by virtue of the maximum principle, $w_1>0$ in $I$ and then it is the first eigenfunction of $-\frac{d^2}{d t^2}$ with Dirichlet boundary conditions in $I$. This in turn implies $\frac{\pi^2}{4} + \alpha=\nu_1(I)= \frac{\pi^2}{4}$, namely $\alpha=0$.
		\end{proof}
	\end{proposizione}
	
	\section{Stability and instability results}\label{section:stability}
	
	In this Section we consider domains, $\Omega_\varphi$, which are hypographs of a positive function, $e^{\varphi}$, in a cylinder and we recall, starting from a positive nondegenerate solution of \eqref{problem_Omega} in $\Omega_\varphi$, how to define an energy functional for small variations of $\varphi$ and so of $\Omega_\varphi$. It is worth to point out that the nondegeneracy of the solution guarantees local uniqueness of the solution under small perturbations of the domain, which is needed to properly define the functional.\\
	Next we introduce the notions of energy-stationary pair and stable/unstable energy-stationary pair under a volume constraint for a couple $(\Omega_\varphi,\varphi)$.\\
	In Subsection \ref{subsection:one-dim sol} we focus on the case when $\Omega_{\varphi_L}$ is a cylinder of height $L$. In particular we first recall some results obtained in \cite{AfonsoPacella} about stability/instability of the couple $(\Omega_{\varphi_L}, u_{p,L})$, where $u_{p,L}$ is the positive one dimensional solution defined in \eqref{upLN}, and then we observe that $(\Omega_{\varphi_L},u_{p,L})$ is a stable/unstable energy pair if and only if so is $(\frac{\Omega_{\varphi_L}}L,u_{p})$, where $u_p:=u_{p,1}$.\\
	At last, we conclude the section, proving Theorem \ref{Teo_stability_p_infty} and Theorem \ref{Teo_stability_p_1}, through the results about the asymptotic analysis of the solution $u_p$, contained in Theorem \ref{Teo_p_to_inf} and Theorem \ref{Teo_p_to_1}.

	\subsection{Energy functional and energy stationary pairs}
	Let $\omega\subset\R^{N-1}$ be a smooth bounded domain and let $\mathcal C_\omega$ be the half cylinder spanned by $\omega$, namely
	\[
	\mathcal C_\omega:=\omega\times(0,+\infty).
	\]
	We denote by $x=(x',x_N)$ the points in $\bar{\mathcal C}_\omega$, where $x'=(x_1,\ldots,x_{N-1})\in\bar\omega$ and $x_N\geq0$.\\
	In $\mathcal C_\omega$ we consider domains whose relative boundaries are cartesian graphs of functions in $C^2(\bar\omega)$. More precisely, for $\varphi\in C^2(\bar\omega)$, we set
	\begin{align*}
		&\Omega_\varphi:=\{(x',x_N)\in\mathcal C_\omega\,:\,x_N<e^{\varphi(x')}\},\\
		&\Gamma_\varphi:=\{(x',x_N)\in\mathcal C_\omega\,:\,x_N=e^{\varphi(x')}\},\\
		&\Gamma_{1,\varphi}:=(\partial \Omega_\varphi\setminus \bar\Gamma_\varphi).
	\end{align*}
	
	We will consider variations of $\Omega_\varphi$ in the class of cartesian graphs of the type $\Omega_{\varphi+tv}$ for $v\in C^2(\bar\omega)$, which amounts to consider a one parameter family of diffeomorphisms $\xi:(-\eta,\eta)\times\bar{\mathcal C}_\omega\to\bar{\mathcal C}_\omega$ of the type
	\[
	\xi(t,x)=(x',e^{tv(x')}x_N),
	\]
	which is a generated by the vector field $V(x)=(0',v(x')x_N)$, where $0'=(0,\ldots,0)\in\R^{N-1}$.\\
	Let $\bar\varphi\in C^2(\omega)$ and let $u_{\Omega_{\bar\varphi}}\in W^{1,\infty}(\Omega_{\bar\varphi})\cap W^{2,2}(\Omega_{\bar\varphi})$ be a positive nondegenerate solution to \eqref{problem_Omega}, with $\Omega=\Omega_{\bar\varphi}$. In \cite[Proposition 2.1]{AfonsoPacella} it has been shown that under such deformations of $\Omega_{\bar\varphi}$ the nondegeneracy of $u_{\Omega_{\bar\varphi}}$ induces a local uniqueness result for solutions of \eqref{problem_Omega} in the deformed domains. Namely, given $v\in C^2(\bar\omega)$, there exists $\delta>0$ such that for any $t\in(-\delta,\delta)$ problem \eqref{problem_Omega} with $\Omega=\Omega_{\bar\varphi+tv}$ admits a unique positive solution $u_{\Omega_{\bar\varphi+tv}}$ in a neighborhood of $u_{\Omega_{\bar\varphi}}\circ\xi(t,\cdot)^{-1}$.
	
	Thus, for any $v\in C^2(\bar\omega)$, the energy functional 
	\[
	T(\Omega_{\bar\varphi+tv})=J(u_{\Omega_{\bar\varphi+tv}})=\frac12\int_{\Omega_{\bar\varphi+tv}}|\nabla u_{\Omega_{\bar\varphi+tv}}(x)|^2\,dx-\frac{1}{p+1}\int_{\Omega_{\bar\varphi+tv}} u_{\Omega_{\bar\varphi+tv}}^{p+1}(x)\,dx
	\]
	is well defined for $t$ sufficiently small.\\
	Ultimately the energy functional $T$ is a functional depending only on functions in $C^2(\bar\omega)$, then, with an abuse of notation, for any $v\in C^2(\bar\omega)$ and for any $t$ sufficiently small we set 
	\[
	T(\bar\varphi+tv):=T(\Omega_{\bar\varphi+tv}).
	\]
	Moreover in \cite[Lemma 4.1]{AfonsoPacella} the first derivative of $T$ at $\bar\varphi$, i.e. for $t=0$, with respect to variations $v \in C^2(\bar\omega)$ has been computed and it takes the following form
	\begin{equation}\label{T'}
		T'(\bar\varphi)[v]=-\frac12\int_\omega \left(\frac{\partial u_{\Omega_{\bar\varphi}}}{\partial \nu}(x',e^{\bar\varphi(x')})\right)^2 v(x') e^{\bar\varphi(x')}dx'.
	\end{equation}
	
	We will be interested in understanding how the energy of a solution $u_{\Omega_{\varphi}}$ behaves with respect to volume-preserving variations of $\Omega_{\varphi}$.\\
	First of all we notice that the volume functional 
	\[
	\mathcal V(\varphi):=\int_\omega e^{\varphi(x')}dx'
	\]
	is of class $C^2$ and $\mathcal V'(\varphi)[v]=\int_\omega e^{\varphi(x')} v(x') dx'$ for any $v\in C^2(\bar\omega)$.\\
	Consequently we define the manifold 
	\[
	M:=\{\varphi\in C^2(\bar\omega)\,:\,\int_{\omega}e^{\varphi (x')}dx'=\int_{\omega}e^{\bar\varphi(x')}dx'\}
	\]
	and we consider the restricted functional
	\[
	I(\varphi):=T_{|M}(\varphi),\quad{\varphi\in M}.
	\]
	Obviously $\bar\varphi\in M$ and the tangent space at $\bar\varphi$ is \[
	T_{\bar\varphi}M=\{v\in C^2(\bar\omega)\,:\,\int_\omega e^{\bar\varphi(x')}v(x')\,dx'=0\}.
	\]
	\begin{definizione}\label{def:critical_energy_pair}
		We say that $(\Omega_{\bar\varphi},u_{\Omega_{\bar\varphi}})$ is an energy-stationary pair under volume constraint if $\bar\varphi$ is a critical point of $I$, namely if $I'(\bar\varphi)[v]=0$ for any $v\in T_{\bar\varphi}M$, or equivalently if there exists $\lambda\in\R$ such that $T'(\bar\varphi)=\lambda \mathcal V'(\bar\varphi)$.
	\end{definizione}
	
	\begin{definizione}\label{def:critical_energy_pair_stable}
		An energy-stationary pair $(\Omega_{\bar\varphi},u_{\Omega_{\bar\varphi}})$ is called stable if 
		\[
		I''(\bar\varphi)[v,v]>0\qquad\forall\:v\in T_{\bar \varphi}M.
		\]
		where $I''$ denotes the second derivative of the restricted functional $I$
	\end{definizione}
	
	\begin{definizione}\label{def:critical_energy_pair_unstable}
		An energy-stationary pair $(\Omega_{\bar\varphi},u_{\Omega_{\bar\varphi}})$ is called unstable if there exists $v\in T_{\bar \varphi}M$ such that
		\[
		I''(\bar\varphi)[v,v]<0.
		\]
	\end{definizione}
	Clearly if $(\Omega_{\bar\varphi},u_{\Omega_{\bar\varphi}})$ is a stable/unstable energy pair, then $\bar\varphi$ is a local minimizer/is not a local minimizer for $T$ under a volume constraint.
	
	\subsection{One-dimensional solutions}\label{subsection:one-dim sol}
	Let $L>0$, let $\varphi_L(x'):=\log L$, $x'\in\bar\omega$, and let us consider 
	\[
	\Omega_{\varphi_L}=\{(x',x_N)\in \mathcal C_\omega\,:\,x_N< L\}.
	\]
	Let $u_{p,L}$ be the positive solution of \eqref{up_problem} and let us extend it to a positive one dimensional solution $u_{p,L}(x',x_N):=u_{p,L}(x_N)$ of \eqref{problem_Omega} in $\Omega_{\varphi_L}$.
	
	By decomposing the spectrum of $L_{u_{p,L}}=-\Delta-pu_{p,L}^{p-1}$, sufficient conditions for nondegeneracy of $u_{p,L}$ have been obtained in \cite[Corollary 4.7]{AfonsoPacella}.
	
	\begin{proposizione}\label{prop:nondegeneracy}
		$u_{p,L}$ is a nondegenerate solution of \eqref{problem_Omega} in $\Omega_{\varphi_L}$ if both the following conditions are satisfied:
		\begin{itemize}
			\item[\emph{i.}] the eigenvalue problem 
			\begin{equation}\label{eigenvalue_problem_alpha_L}
				\begin{cases}
					-z'' - pu_{p,L}^{p-1}z=\alpha z \quad \mbox{ in } (0,L)\\
					z'(0)=z(L)=0
				\end{cases}
			\end{equation}
			does not admit zero as an eigenvalue;
			\item[\emph{ii.}] $\lambda_1(\omega)>-\alpha_{1,L}(p)$, where $\alpha_{1,L}(p)$ is the first eigenvalue of \eqref{eigenvalue_problem_alpha_L}.
		\end{itemize}   
	\end{proposizione}
	\begin{osservazione}\label{rem:nondegeneracy}
		Condition {i.} is always satisfied by virtue of Proposition \ref{prop:nondeg}, so if $\lambda_1(\omega)>-\alpha_{1,L}(p)$ the one-dimensional solution $u_{p,L}$ to \eqref{problem_Omega} in $\Omega_{\varphi_L}$ is nondegenerate.    
	\end{osservazione}
	
	As a consequence, if we define
	\[
	M=\{\varphi\in C^2(\bar\omega)\,:\,\int_\omega e^{\varphi(x')}dx'=\int_\omega e^{\varphi_L(x')}dx'=L|\omega|\},
	\]
	then $(\Omega_{\varphi_L}, u_{p,L})$ is an energy-stationary pair under volume constraint, being, by \eqref{T'},
	\[
	T'(\varphi_L)[v]=-\frac{L}{2}(u'_{p,L})^2\int_\omega v(x')\,dx'=0
	\]
	for any $v\in T_{\varphi_L}M=\{v\in C^2(\bar\omega)\,:\,\int_\omega v(x')\,dx'=0\}$.
	
	Next we recall a result obtained in \cite{AfonsoPacella}, that will be crucial in order to prove Theorem \ref{Teo_stability_p_infty} and Theorem \ref{Teo_stability_p_1}. 
	\begin{teorema}\label{teorema}
		Let $\omega\subset\R^{N-1}$ be a smooth bounded domain.\\
		Let $u_{p,L}$ be the positive one-dimensional solution to \eqref{problem_Omega} in $\Omega_{\varphi_L}$. Let $\lambda_1(\omega)$ be the first nontrivial Neumann eigenvalue of $-\Delta_{\R^{N-1}}$ in $\omega$, let $\alpha_{1,L}(p)$ be the first eigenvalue of \eqref{L_up} with Dirichlet boundary conditions in $(-L,L)$ and assume that $\lambda_1(\omega)>-\alpha_{1,L}(p)$.\\
		Let $h_{p,L}$ be the solution to 
		\begin{equation}\label{hpLequation}
			\begin{cases}
				h''_{p,L}=(\lambda_1(\omega)-pu_{p,L}^{p-1})h_{p,L} \quad \mbox{ in } (0,L)\\
				h_{p,L}(L)=-Lu'_{p,L}(L)\\
				h'_{p,L}(0)=0.
			\end{cases}
		\end{equation}
		Then:\\
		\emph{$i)$} $(\Omega_{\varphi_L}, u_{p,L})$ is an unstable energy stationary pair if and only if $h'_{p,L}(L)<0$;\\
		\emph{$ii)$} $(\Omega_{\varphi_L}, u_{p,L})$ is a stable energy stationary pair if and only if $h'_{p,L}(L)>0$.   
	\end{teorema}
	In fact in \cite{AfonsoPacella} this result has been obtained in the case $L=1$, but it can be easily extended to the general case $L>0$, and more importantly in \cite{AfonsoPacella} these conditions are stated just as sufficient ones, but actually they hold as \virg{if and only if}, as it is clear from the proofs of \cite[Theorem 4.11]{AfonsoPacella} and \cite[Theorem 1.5]{AfonsoPacella}.
	
	It is worth pointing out that assumption $\lambda_1(\omega)>-\alpha_{1,L}(p)$ is not only a sufficient condition to guarantee the nondegeneracy of the one-dimensional solution $u_{p,L}$ of \eqref{problem_Omega} in $\Omega_{\varphi_L}$ (see Remark \ref{rem:nondegeneracy}), but also implies the positivity of the function $h_{p,L}$, which is a key property in the proof of Theorem \ref{teorema}. We state the positivity result, obtained in  \cite[Proposition 4.10]{AfonsoPacella} via maximum principle, in a suitable form and for $L=1$. 
	\begin{proposizione}\label{prop:hpL>0}
		Let $\lambda>0$ and let assume that $\lambda>-\alpha_{1}(p)$, where $\alpha_{1}(p)$ is the first eigenvalue of \eqref{L_up} with Dirichlet boundary conditions in $(0,1)$\\
		Then the solution $\tilde h_{p}$ to 
		\begin{equation}\label{hpLequation}
			\begin{cases}
				\tilde h''_{p}=(\lambda-pu_{p}^{p-1}) \tilde h_{p} \quad \mbox{ in } (0,1)\\
				\tilde h_{p}(1)=-u'_{p}(1)\\
				\tilde h'_{p}(0)=0.
			\end{cases}
		\end{equation}
		is positive in $[0,1]$.
	\end{proposizione}

	\begin{osservazione}\label{rem:equivalence}
		Since $u_{p,L}(y)=\left(\frac1L\right)^{\frac{2}{p-1}} u_p(\frac{y}{L})$,\;\;$y\in(0,L)$, and we have studied the asymptotic behaviour {of the solution $u_p$ of \eqref{problema_in_I}} in (0,1),
		it is convenient for our purposes to define in $(0,1)$ a rescaling $h_p$ of $h_{p,L}$, namely
		\[
		h_p(t):=L^{\frac{2}{p-1}}h_{p,L}(Lt),\quad t\in(0,1),
		\]
		which solves 
		\begin{equation}\label{hpequation}
			\begin{cases}
				h''_{p}=(L^2\lambda_1(\omega)-pu_{p}^{p-1})h_{p} \quad \mbox{ in } (0,1)\\
				h_{p}(1)=-u'_{p}(1)\\
				h'_{p}(0)=0.
			\end{cases}
		\end{equation}
		Being 
		\begin{equation}\label{rescalings}
			\lambda_1\left(\frac{\omega}{L}\right)=L^2\lambda_1(\omega)\qquad\text{and}\qquad\alpha_{1}(p)=L^2\alpha_{1,L}(p),
		\end{equation}
		where $\alpha_1(p):=\alpha_{1,1}(p)$, then
		\begin{equation}\label{iffs}
			\lambda_1(\omega)+\alpha_{1,L}(p)>0\quad \Leftrightarrow\quad \lambda_1\left(\frac{\omega}{L}\right)+\alpha_1(p)>0
			\quad \Leftrightarrow\quad L^2 \lambda_1(\omega)+\alpha_1(p)>0.
		\end{equation}
		Hence, if we assume $\lambda_1(\omega)+\alpha_{1,L}(p)>0$, $h_p$ is the unique solution to \eqref{hpequation} and $h'_{p,L}(L)\gtrless0$ if and only if $h'_p(1)\gtrless 0$. Therefore the sign of $h'_p(1)$ is what determines the stability/instability, by \emph{i)}, \emph{ii)} of Theorem \ref{teorema}.
		Furthermore by Theorem \ref{teorema}, \eqref{hpequation}, \eqref{rescalings} and \eqref{iffs} we have that $(\Omega_{\varphi_L}, u_{p,L})$ is a stable/unstable energy stationary pair if and only if $(\frac{\Omega_{\varphi_L}}{L},u_p)$ is a stable/unstable energy stationary pair.
	\end{osservazione}

	\subsection{Proof of Theorem \ref{Teo_stability_p_infty}} 
	By Theorem \ref{teorema} and Remark \ref{rem:equivalence} it is sufficient to prove that, for $p$ sufficiently large and $L^2\lambda_1(\omega)>\gamma^2p^2$, we have $h'_p(1)>0$, where $h_p$ is the solution to \eqref{hpequation}.\\
	We start by noticing that by Proposition \ref{prop:hpL>0}, applied with $\lambda=L^2\lambda_1(\omega)$, the function $h_p$ is positive in $[0,1]$. Then let us consider the even extension of $h_p$ to $\bar I=[-1,1]$, that, with an abuse of notation, we will still denote as $h_p$.\\
	Next in order to simplify the notation we set
	\[
	\lambda:=L^2\lambda_1(\omega).
	\]
	If $\lambda \geq p\|u_p\|_\infty^{p-1}$ then, for any $p>1$, by $\eqref{hpLequation}$ we have $h''_p>0$ in $(0,1)$, which implies $h'_p(1)>h'_p(0)=0$. This case was already covered by \cite{AfonsoPacella}.\\

	Whereas we claim that for $p$ sufficiently large, if $\lambda\in [{\gamma^2 p^2},p\|u_p\|_\infty^{p-1})$
	\begin{equation}\label{claim:mp}
		\exists\: m_p\in(0,1)\qquad\text{such that}\quad h'_p<0 \mbox{ in } (0,m_p) \quad h'_p>0 \mbox{ in } (m_p,1].    
	\end{equation}
	We notice that the interval ${[\gamma^2 p^2},p\|u_p\|_\infty^{p-1})$ is not empty by {Proposition \ref{prop:Lp_estimate}}.\\
	Clearly, once the claim is proved the thesis follows.\\
	Suppose by contradiction that there exists $p_n\to+\infty$ and $\lambda_n\in {[\gamma^2 p_n^2},p_n\|u_{p_n}\|_\infty^{p_n-1})$ such that $\|h_{p_n}\|_\infty=h_{p_n}(0)$, then
	\begin{equation*}
		\frac{\gamma^2p_n^2}{ p_n\|u_{p_n}\|_\infty^{p_n-1}}\leq\mu_{p_n}^2{\lambda_n}<1,
	\end{equation*}
	where $\mu_{p_n}=(p_n \|u_{p_n}\|_\infty^{p_n-1})^{-\frac12}$.\\
	So, by $(3)$ of Theorem \ref{Teo_p_to_inf}, up to a subsequence, ${\lambda_n}\mu_{p_n}^2\rightarrow {\eta\in[2\gamma^2,1]\subset(\frac12,1]}$, as $n\to +\infty$, and rescaling the function by 
	\begin{equation*}
		\tilde{h}_{p_n}(y)=\frac{h_{p_n}(\mu_{p_n} y)}{h_{p_n}(0)}, \qquad y\in \left[ -\frac{1}{\mu_{p_n}},\frac{1}{\mu_{p_n}}\right],
	\end{equation*}
	we get that $\Vert \tilde{h}_{p_n}\Vert_{\infty}=\tilde{h}_{p_n}(0)=1$ and that $\tilde{h}_{p_n}$ solves the equation
	{
		\begin{equation*}
			\begin{cases}
				\tilde{h}''_{p_n}(y)=\left( \mu_{p_n}^2\lambda_n-\left(1+\frac{\tilde{u}_{p_n}(y)}{p_n}\right)^{p_n-1}\right)\tilde{h}_{p_n}(y) \qquad y\in \left( -\frac{1}{\mu_{p_n}},\frac{1}{\mu_{p_n}}\right)\\
				\tilde{h}_{p_n}(0)=1\\
				\tilde{h}'_{p_n}(0)=0\\
				\tilde{h}_{p_n}>0. 
			\end{cases}
		\end{equation*}
	}
	Moreover, by Lemma \ref{lemma:decay}, we have 
	\begin{align*}
		|\tilde{h}_{p_n}(y)|&\leq 1,\\
		|\tilde{h}'_{p_n}(y)|&\leq \int_{0}^{y}\left|\mu_{p_n}\lambda_n-\left(1+\frac{\tilde{u}_{p_n}(t)}{p_n}\right)^{p_n-1}\right|\tilde{h}_{p_n}(t)\,dt\\
		&\leq \int_{0}^{y}\left(\mu_{p_n}\lambda_n+e^{\frac{p_n-1}{p_n}\tilde u_{p_n}(t)}\right)\tilde{h}_{p_n}(t)\,dt\\
		&\leq \int_{0}^{y} 1 +e^{\frac{p_n-1}{p_n}(C_1W(t)+C_2)}\,dt \leq C 
	\end{align*}
	in every compact subset of $\left[ -\frac{1}{\mu_{p_n}},\frac{1}{\mu_{p_n}}\right]$, where $W$ is defined in \eqref{W}. Of course, using the equation and repeating the previous calculations, it is immediate to see that $\tilde{h}''_{p_n}$ is bounded too. So we have
	\begin{equation*}
		\tilde{h}_{p_n} \rightarrow \tilde{H} \quad \mbox{ in } C^1_{loc}(\mathbb{R})
	\end{equation*}
	where $\tilde{H}$ solves
	\begin{equation*}
		\begin{cases}
			-\tilde{H}''=e^W\tilde{H}-{\eta} \tilde{H}\qquad\text{in $\R$}\\
			\tilde{H}(0)=1, \tilde{H}'(0)=0\\
			\|\tilde H\|_\infty=1.
		\end{cases}
	\end{equation*}

	{
		Reasoning exactly as in the last part of the proof of Theorem $\ref{thm:alpha1p}$ we conclude that $\eta=\frac{1}{2}$, which is a contradiction against $\eta \in (\frac12,1]$.}
	\hfill \qed\\
	
	\subsection{Proof of Theorem \ref{Teo_stability_p_1}}
	First of all by \eqref{rescalings} and $(4)$ Theorem \ref{Teo_p_to_1}, as $p\to1$,
	\begin{equation*}
		\lambda_1(\omega)+\alpha_{1,L}(p)=\lambda_1(\omega)+\frac{\alpha_{1}(p)}{L^2}\rightarrow \lambda_1(\omega) >0.
	\end{equation*} 
	As a consequence, for any $L>0$ and $p$ sufficiently close to $1$, $u_{p,L}$ is a nondegenerate solution to \eqref{problem_Omega} in $\Omega_{\varphi_L}$ by Remark \ref{rem:nondegeneracy}.
	Next, in order to simplify the notation we set
	\[
	\lambda:=L^2\lambda_1(\omega).
	\]
	To prove the stability/instability, by Remark \ref{rem:equivalence} it is sufficient to prove that, for $\lambda\gtrless\frac{\pi^2}{4}$ and $p$ sufficiently close to $1$, $h'_p(1)\gtrless0$, where $h_p$ is the solution to \eqref{hpequation}.\\
	We start by \eqref{iffs} noticing that, by Proposition \ref{prop:hpL>0}, $h_p$ is positive in $[0,1]$. Then let us consider the even extension of $h_p$ to $I=(-1,1)$, that, with an abuse of notation, we will still denote as $h_p$.\\
	
	If $\lambda > \frac{\pi^2}{4}$, then by \eqref{convergence_from_above}, for $p$ sufficiently close to $1$, $\lambda> pu_p^{p-1}$ in $[0,1]$ and so, by $\eqref{hpLequation}$, for $p$ sufficiently close to $1$,  
	$h''_p(t)>0$ for all $t \in [0,1]$, which implies $h'_p(1)>h'_p(0)=0$.\\
	Whereas if $\lambda < \frac{\pi^2}{4}$, by \eqref{convergence_from_above} we deduce the existence of $t_p\in(0,1)$ such that 
	\begin{equation}\label{inflection_point}
		\lambda=pu_p^{p-1}(t_p).
	\end{equation}
	So $t_p$ is an inflection point and
	\begin{equation*}
		h''_p<0 \quad \mbox{ in }(0,t_p), \quad h''_p>0 \quad \mbox{ in } (t_p,1).
	\end{equation*}
	
	Let us show that 
	\begin{equation}
		\label{tp_to_1}
		t_p\to1\qquad\text{as $p\to1$.}
	\end{equation}
	If this is not the case, then there exists $p_n \to 1$, as $n\to+\infty$, and $\varepsilon>0$ such that $t_{p_n}\in(0,1-\varepsilon]$ for any $n$ sufficiently large. Then by $(3)$ of Theorem \ref{Teo_p_to_1}
	\[
	\lambda=p_n u_{p_n}^{p_n-1}(t_n)\to\frac{\pi^2}{4}\qquad\text{as $n\to\infty$},
	\]
	which is a contradiction against $\lambda<\frac{\pi^2}{4}$.\\
	Next we consider $\bar{h}_p=\frac{h_p}{\|u_p\|_\infty}$, which solves
	\begin{equation*}
		\begin{cases}
			\bar{h}''_p=(\lambda-pu_p^{p-1})\bar{h}_p \qquad \mbox{ in } I\\
			\bar{h}_p(\pm1)=-\bar{u}'_p(1)=\frac{\pi}{2}(1+o_p(1))\\
			\bar{h}_p'(0)=0,
		\end{cases}
	\end{equation*}
	where $\bar{u}_p=\frac{u_p}{\|u_p\|_\infty}$ and $-\bar{u}'_p(1)=\frac{\pi}{2}(1+o_p(1))$ by $(2)$ of Theorem \ref{Teo_p_to_1}.\\
	
	We are now in position to prove the following {claim}:
	\begin{equation}\label{claim:hp0}
		\bar{h}_p(0) \nrightarrow 0 \quad \mbox{ when } p\to 1.
	\end{equation}
	Indeed if by contradiction we have
	\begin{equation*}
		\bar{h}_p(0) \rightarrow 0 
	\end{equation*}
	then 
	\begin{equation*}
		0<\bar{h}_p(t_p)<\bar{h}_p(0)\rightarrow 0 \quad \mbox{ because } h'_p<0  \quad \mbox{ in } (0,t_p),
	\end{equation*}
	which means $\bar{h}_p(t_p) \rightarrow 0$. Moreover, by \eqref{convergence_from_above} and being $\|\bar{h}_p\|_{\infty}=\bar{h}_p(1)=\frac{\pi}{2}(1+o_p(1))$, we get,
	\begin{equation*}
		|\bar{h}'_p(t)| =\left| \int_{0}^{t} (\lambda-pu_p^{p-1}(s))\bar{h}_p(s)\,ds \right| 
		\leq \int_{0}^{1} (\lambda+pu_p^{p-1}(s))\|\bar{h}_p\|_{\infty}\,ds\leq C.
	\end{equation*}
	In turn this gives a contradiction, indeed
	\begin{equation*}
		\frac{\pi}{2}(1+o_p(1))=\bar{h}_p(1)=\bar{h}_p(t_p)+\int_{t_p}^{1}\bar{h}'_p(s)\,ds \leq o_p(1)+C(1-t_p)\overset{\eqref{tp_to_1}}{=}o_p(1),
	\end{equation*}
	and this proves \eqref{claim:hp0}.
	
	Let $\tilde{h}_p=\frac{\bar{h}_p}{\bar{h}_p(0)}$ then 
	\begin{equation*}
		\tilde{h}_p(0)=1, \quad
		\|\tilde{h}_p\|_{\infty}= \max \biggl\{ 1, \frac{\bar{h}_p(1)}{\bar{h}_p(0)} \biggr\} \leq C.
	\end{equation*}
	Let $\eta \in (0,1)$ be fixed, then by $(3)$ of Theorem \ref{Teo_p_to_1} we have $pu_p^{p-1}\rightrightarrows\frac{\pi^2}{4}$ in $[-\eta,\eta]$. Furthermore
	\begin{equation*}
		\begin{cases}
			\tilde{h}''_p=(\lambda-pu_p^{p-1})\tilde{h}_p \quad \mbox{ in }[-\eta,\eta]\\
			\tilde{h}_p(0)=1,\\
			\tilde{h}_p'(0)=0,
		\end{cases}
	\end{equation*}
	then both $\tilde{h}_p$ and its derivatives are bounded in $[-\eta,\eta]$, indeed
	\begin{align}
		&\|\tilde{h}_p\|_{\infty} \leq C,\label{last}\\
		&| \tilde{h}''_p| \leq (\lambda+\frac{\pi^2}4+o_p(1))C\leq C,\nonumber\\
		&|\tilde{h}'_p(t)|\leq \left|\int_{0}^{t}|\tilde{h''}_p(\tau)|\,d \tau\right|\leq C.\nonumber
	\end{align}
	Therefore by Ascoli-Arzelà theorem $\tilde{h}_p \rightrightarrows H_\lambda$ in $C^1([-\eta,\eta])$ which solves
	\begin{equation*}
		\begin{cases}
			H''_\lambda=\left(\lambda-\frac{\pi^2}{4}\right) H_\lambda \quad \mbox{ in } (-\eta,\eta)\\
			H_\lambda(0)=1\\
			H_\lambda'(0)=0.
		\end{cases}
	\end{equation*}
	The solution of this equation is easily computable and it is
	\begin{equation*}
		H_\lambda(t)=\cos\left( \sqrt{\frac{\pi^2}{4}-\lambda}\:t\right), \quad t \in [-\eta,\eta].
	\end{equation*}
	Hence, since $t_p>\eta$ for $p$ sufficiently close to $1$,  we have by $(3)$ of Theorem \ref{Teo_p_to_1} and \eqref{last}
	\begin{align*}
		\tilde{h}'_p(1)&<\int_{0}^{\eta} (\lambda-pu_p^{p-1}(t))\tilde{h}_p(t)\,dt + \int_{t_p}^{1} (\lambda-pu_p^{p-1}(t))\tilde{h}_p(t)\,dt\\
		&<\int_{0}^{\eta} (\lambda-\frac{\pi^2}{4})\cos\left( \sqrt{\frac{\pi^2}{4}-\lambda}\:t\right)\,dt+o_p(1) + \int_{t_p}^{1} \lambda\|\tilde{h}_p\|_{\infty}\,dt\\
		&=\int_{0}^{\eta} (\lambda-\frac{\pi^2}{4})\cos\left( \sqrt{\frac{\pi^2}{4}-\lambda}\,t\right)\,dt+o_p(1) +(1-t_p)\lambda C\\
		&=- \sqrt{\frac{\pi^2}{4}-\lambda}\sin\left( \sqrt{\frac{\pi^2}{4}-\lambda}\:\eta\right)+o_p(1)<0 \quad \mbox{ as } p\to 1.
	\end{align*}
	And this is enough to conclude, being
	\begin{equation*}
		h'_p(1)=h_p(0)\tilde{h}'_p(1)<0.
	\end{equation*}
	\hfill \qed\\

	\

\end{document}